\documentclass{article}
\usepackage[utf8]{inputenc}

\usepackage{authblk}
\usepackage{hyperref}
\usepackage{comment}
    
\title{Interpolating between BSDEs and PINNs: \\deep learning for elliptic and parabolic boundary value problems}
\author[1]{Nikolas N\"usken}
\author[2,3,4]{Lorenz Richter}
\date{\today}

\affil[1]{Department of Mathematics, King's College London, UK, \href{mailto:nikolas.nusken@kcl.ac.uk}{nikolas.nusken@kcl.ac.uk}}
\affil[2]{dida Datenschmiede GmbH, 10827 Berlin, Germany, \href{mailto:lorenz.richter@dida.do}{lorenz.richter@dida.do}}
\affil[3]{Institute of Mathematics, Freie Universit\"at Berlin, 14195 Berlin, Germany}
\affil[4]{Zuse Institute Berlin, 14195 Berlin, Germany}

\usepackage{mathrsfs, amssymb}
\usepackage{amsthm}
\usepackage[intlimits]{empheq}
\usepackage{bbm}     
\usepackage{amsmath}
\usepackage{xcolor}
\usepackage[scientific-notation=true]{siunitx}
\usepackage{float}
\usepackage{dsfont}	
\usepackage[nottoc,numbib]{tocbibind}
\usepackage{bbm}
\usepackage{subcaption}
\usepackage[inner=2cm,outer=2cm,top=2cm,bottom=2cm]{geometry}
\usepackage[ruled,vlined]{algorithm2e}
\usepackage{cleveref}
\usepackage{caption, subcaption}
\usepackage{pifont}
\usepackage{enumerate}   
\usepackage{todonotes}
\usepackage{tabularx}
\usepackage[maxcitenames=8,backend=bibtex,style=trad-abbrv,sorting=nyt,sortcites]{biblatex}
\newcommand{\xmark}{\ding{55}}

\let\P\relax\newcommand{\P}{\mathbb{P}}

\DeclareMathOperator{\E}{\mathbb{E}}
\DeclareMathOperator{\R}{\mathbb{R}}

\theoremstyle{plain}
\newtheorem{theorem}{Theorem}[section]

\newtheorem{proposition}[theorem]{Proposition}
\newtheorem{assumption}{Assumption}

\theoremstyle{definition}
\newtheorem{definition}[theorem]{Definition}

\theoremstyle{remark}
\newtheorem{remark}[theorem]{Remark}

\usepackage[T1]{fontenc}
\usepackage{lipsum}
\usepackage{tocloft}
\begin{document}
    
\maketitle

\begin{abstract}
Solving high-dimensional partial differential equations is a recurrent challenge in economics, science and engineering. In recent years, a great number of computational approaches have been developed, most of them relying on a combination of Monte Carlo sampling and deep learning based approximation. For elliptic and parabolic problems, existing methods can broadly be classified into those resting on reformulations in terms of \emph{backward stochastic differential equations} (BSDEs) and those aiming to minimize a regression-type $L^2$-error (\emph{physics-informed neural networks}, PINNs). In this paper, we review the literature and suggest a methodology based on the novel \emph{diffusion loss} that interpolates between BSDEs and PINNs. Our contribution opens the door towards a unified understanding of numerical approaches for high-dimensional PDEs, as well as for implementations that combine the strengths of BSDEs and PINNs. The diffusion loss furthermore bears close similarities to \emph{(least squares) temporal difference} objectives found in reinforcement learning. We also discuss eigenvalue problems and perform extensive numerical studies, including calculations of the ground state for nonlinear Schr\"odinger operators and committor functions relevant in molecular dynamics. 
\end{abstract}

\section{Introduction}

In this article, we consider approaches towards solving high-dimensional partial differential equations (PDEs) that are based on minimizing appropriate loss functions in the spirit of machine learning. For example, we aim at identifying approximate solutions to nonlinear parabolic boundary value problems of the form
\begin{subequations}
\label{eq: general boundary value problem}
\begin{align}
    \label{eq: general PDE}
    (\partial_t + L) V(x, t) + h(x, t, V(x, t), \sigma^\top\nabla V(x, t)) &= 0,\qquad  & (x, t) \in \Omega \times [0, T), \\
    \label{eq: time boundary condition}
    V(x, T) & = f(x), \qquad  & x \in \Omega,\\
    \label{eq: space boundary condition}
    V(x, t) & = g(x, t), \qquad  & (x, t) \in \partial \Omega \times [0, T],
\end{align}
\end{subequations}
on a spatial domain $\Omega \subset \R^d$ and time interval $[0,T]$, where $h \in C(\R^d \times [0, T] \times \R \times \R^{d}, \R)$ specifies the nonlinearity, and $f \in C(\R^d, \R)$ as well as $g \in C(\R^d \times [0, T], \R)$ are given functions defining the terminal and boundary conditions. Moreover,
\begin{equation}
\label{eq: infinitesimal generator variational}
    L = \frac{1}{2} \sum_{i,j=1}^d (\sigma \sigma^\top)_{ij}(x,t) \partial_{x_i} \partial_{x_j} + \sum_{i=1}^d b_i(x,t) \partial_{x_i}
\end{equation}
is an elliptic differential operator including the coefficient functions $b \in C(\R^d \times [0, T], \R^d)$ and $\sigma \in C(\R^d \times [0, T], \R^{d \times d})$, with $\sigma$ assumed to be non-degenerate. We will later make use of the fact that $L$ is the infinitesimal generator of the diffusion process defined by the stochastic differential equation (SDE)
\begin{equation}
\label{eq: uncontrolled SDE}
    \mathrm dX_s = b(X_s, s) \,\mathrm d s + \sigma(X_s, s) \,\mathrm d W_s,
\end{equation}
where $W_s$ is a standard $d$-dimensional Brownian motion. Our approach exploiting the connection between \eqref{eq: general boundary value problem} and \eqref{eq: uncontrolled SDE} extends almost effortlessly to a wide range of nonlinear elliptic PDEs (see \Cref{sec: special PDEs}) and eigenvalue problems (see \Cref{sec: eigenvalue problems}).
We also note that the terminal condition \eqref{eq: time boundary condition} can be replaced by an initial condition without loss of generality, using a time reversal transformation. Similarly, we may replace the Dirichlet boundary condition \eqref{eq: space boundary condition} by its Neumann counterpart, that is, constraining the normal derivative $\partial_{\vec{n}} V$ on $\partial \Omega \times [0,T]$ in \eqref{eq: space boundary condition}.
\par\bigskip

A notorious challenge that appears in the numerical treatment of PDEs is the \emph{curse of dimensionality}, suggesting that the computational complexity increases exponentially in the dimension of the state space. In recent years, however, multiple numerical \cite{hure2020deep,raissi2017physics,weinan2017deep} as well as theoretical studies \cite{grohs2019deep,jentzen2018proof} have indicated that a combination of Monte Carlo methods and neural networks offers a promising way to overcome this problem. This paper centers around two strategies that allow for solving quite general nonlinear PDEs:
\begin{itemize}
    \item 
    \textbf{PINNs (physics-informed neural networks)} \cite{raissi2017physics}, also known under the name DGM \emph{(deep Galerkin method)} \cite{sirignano2018dgm} directly minimize the misfit between the left-hand sides and right-hand sides of \eqref{eq: general boundary value problem}, evaluated at appropriately chosen (random) points in the space-time domain $\Omega \times [0,T]$ and its boundary.
    \item
    \textbf{deep BSDEs (backward stochastic differential equations)} \cite{weinan2017deep} rely on a reformulation of \eqref{eq: general boundary value problem} in terms of stochastic forward-backward dynamics, explicitly making use of the connection between the PDE \eqref{eq: general boundary value problem} and the SDE \eqref{eq: uncontrolled SDE}. Deep BSDEs minimize the misfit in the terminal condition associated to the backward SDE. 
\end{itemize}
We review those approaches (see Section \ref{sec: PINN and BSDE loss}) and -- motivated by It\^o's formula -- introduce a novel optimization objective, called the \emph{diffusion loss} $\mathcal{L}_{\mathrm{diffusion}}^\mathfrak{t}$ (see Section \ref{sec: diffusion loss}). Importantly, our construction depends on the auxiliary time parameter $\mathfrak{t} \in (0, \infty)$, allowing us to recover PINNs in the limit $\mathfrak{t} \rightarrow 0$ and deep BSDEs in the limit $\mathfrak{t} \rightarrow \infty$:
$$
\text{PINNs} \xleftarrow{\quad\mathfrak{t} \rightarrow 0 \quad}\mathcal{L}_{\mathrm{diffusion}}^\mathfrak{t} \xrightarrow{\quad \mathfrak{t} \rightarrow \infty \quad } \text{deep BSDEs}
$$
As will become clear below, PINNs may therefore be thought of as accumulating derivative information on the PDE solution locally (in time), whereas deep BSDEs constitute global approximation schemes (relying on entire trajectories). Moreover, the diffusion loss is closely related to \emph{(least squares) temporal difference learning} (see \cite[Chapter 3]{szepesvari2010algorithms} and Remark \ref{rem:reinforcement} below).\footnote{We would like to thank an anonymous referee for pointing out this connection.}
\par\bigskip

The diffusion loss provides an interpolation between seemingly quite distinct methods. Besides this theoretical insight, we show experimentally that an appropriate choice of $\mathfrak{t}$ can lead to a computationally favorable blending of PINNs and deep BSDEs, combining the advantages of both methods. In particular, the diffusion loss (with $\mathfrak{t}$ chosen to be of moderate size) appears to perform well in scenarios where the domain $\Omega$ has a (possibly complex) boundary and/or the PDE under consideration contains a large number of second order derivatives (for instance, when $\sigma$ is not sparse). We will discuss the trade-offs involved in Section \ref{sec: from losses to algorithms}.
\par\bigskip

\par\bigskip

\textbf{The curse of dimensionality, BSDEs vs. PINNs.}
Traditional numerical methods for solving PDEs (such as finite difference and  finite volume methods, see \cite{ames2014numerical}) usually require a discretization of the domain $\Omega$, incurring a computational  cost that for a prescribed accuracy is expected to grow linearly in the number of grid points, and hence exponentially in the number of dimensions. The recently developed  approaches towards beating this curse of dimensionality  -- BSDEs and PINNs alike\footnote{We note in passing that multilevel Picard approximations \cite{becker2020numerical,hutzenthaler2020overcoming,hutzenthaler2016multilevel} represent another interesting and fairly different class of methods that however are beyond the scope of this paper.} --  replace the deterministic mesh by Monte Carlo sampling, in principle promising dimension-independent convergence rates \cite{e2020algorithms}. Typically, the approximating function class is comprised of neural networks (expectantly providing adaptivity to low-dimensional latent structures \cite{bach2017breaking}), although also tensor trains have shown to perform well \cite{richter2021solving}. In contrast to PINNs, methods based on BSDEs make essential use of the underlying diffusion \eqref{eq: uncontrolled SDE} to generate the training data. According to our preliminary numerical experiments, it is not clear whether this use of structure indeed translates into computational benefits; we believe that additional research into this comparison is needed and likely to contribute both conceptually and practically to further advances in the field. The diffusion loss considered in this paper may be a first step in this direction, for a direct comparison between BSDEs and PINNs we refer to
\Cref{tab: loss comparisons data generation} and 
\Cref{tab: loss comparisons challenges} below.

\par\bigskip

\textbf{Previous works.} Attempts to approximate PDE solutions by combining Monte Carlo methods with neural networks date back to the 1990s, where some variants of residual minimizations have been suggested \cite{lagaris1998artificial, lagaris2000neural, uchiyama1993solving, lee1990neural, jianyu2003numerical}. Recently, this idea gained popularity under the names \textit{physics-informed neural networks} (PINNs, \cite{raissi2017physics, raissi2019physics}) and \textit{deep Galerkin methods} (DGMs, \cite{sirignano2018dgm}). Let us further refer to a comparable approach that has been suggested in \cite{berg2018unified}, and, for the special case of solving the dynamic Schr\"odinger equation, to \cite{carleo2017solving}. For theoretical analyses we shall for instance mention \cite{mishra2020estimates}, which provides upper bounds on the generalization error, \cite{de2021error}, which states an error analysis for PINNs in the approximation of Kolmogorov-type PDEs, \cite{muller2021notes}, which investigates convergences depending on whether using exact or penalized boundary terms, and \cite{wang2022and, wang2020eigenvector}, which study convergence properties through the lens of neural tangent kernels. Some further numerical experiments have been conducted in \cite{dockhorn2019discussion, magill2018neural} and multiple algorithmic improvements have been suggested, e.g. in \cite{wang2021understanding}, which balances gradients during model training, \cite{jagtap2020adaptive}, which considers adaptive activation functions to accelerate convergence, as well as \cite{van2020optimally}, which investigates efficient weight-tuning.
\par\bigskip
BSDEs have first been introduced in the 1970s \cite{bismut1973conjugate} and eventually studied more systematically in the 1990s \cite{pardoux1990adapted}. For a comprehensive introduction elaborating on their connections to both elliptic and parabolic PDEs we refer for instance to \cite{pardoux1998backward}. Numerical attempts to exploit this connection aiming for approximations of PDE solutions have first been approached by backward-in-time iterations, originally relying on a set of basis functions and addressing parabolic equations on unbounded domains \cite{pardoux1998backward, bouchard2004discrete, gobet2005regression}. Those methods have been considered and further developed with neural networks in \cite{hure2020deep, beck2019deep} and endowed with tensor trains in \cite{richter2021solving}. A variational formulation termed \textit{deep BSDE} has been first introduced in \cite{weinan2017deep, han2018solving}, aiming at PDE solutions at a single point, with some variants following e.g. in \cite{raissi2018forward, nusken2021solving, richter2021solvingPDEs}. For an analysis of the approximation error of the deep BSDE method we refer to \cite{han2020convergence} and for a fully nonlinear version based on second-order BSDEs to \cite{beck2019machine}. An extension to elliptic PDEs on bounded domains has been suggested in \cite{kremsner2020deep}.
\par\bigskip
The diffusion loss as well as other methods that aim at solving high-dimensional PDEs can be related to well-known approaches in reinforcement learning and optimal control. In particular, the key definition \eqref{eq: definition diffusion loss int} below can be seen as a generalization of equation (2.27) in \cite{zhou2021actor}, derived in the context of (least squares) temporal difference learning \cite[Chapter 3]{szepesvari2010algorithms}. In contrast, the algorithms suggested in \cite{han2020derivative} and  \cite{martin2022solving} are in the spirit of Monte Carlo reinforcement techniques \cite[Section 3.1.2]{szepesvari2010algorithms}, relying on an averaged rather than a pathwise construction principle for optimization objectives. We would also like to point out that introducing the artificial time parameter $\mathfrak{t}$ is conceptually similar to the idea in \cite{han2020solving} where the authors elevate an eigenvalue problem of elliptic type to a parabolic one using a related augmentation.

\par\bigskip
There are additional works aiming to solve linear PDEs specifically, mainly by exploiting the Feynman-Kac theorem and mostly considering parabolic equations \cite{beck2018solving, berner2018analysis, richter2022robust}. For the special case of the Poisson equation, \cite{grohs2020deep} considers an elliptic equation on a bounded domain. Linear PDEs often admit a variational formulation that suggests to minimize a certain energy functional -- this connection has been used in \cite{weinan2018deep, khoo2019solving, lu2021priori} and we refer to \cite{muller2019deep} for some further analysis. Similar minimization strategies can be used when considering eigenvalue problems, where we mention \cite{zhang2021solving} in the context of metastable diffusion processes. We also refer to \cite{lagaris1997artificial, jin2020unsupervised, pfau2020ab, hermann2020deep} for similar problems in quantum mechanics that often rest on particular neural network architectures. Nonlinear elliptic eigenvalue problems are addressed in \cite{pfau2020ab, han2020solving} by exploiting a connection to a parabolic PDE and a fixed point problem.
\par\bigskip
For rigorous results towards the capability of neural networks to overcome the curse of dimensionality we refer to \cite{jentzen2018proof, grohs2018proof, hutzenthaler2020proof}, each analyzing certain special PDE cases. Adding to the methods referred to above, let us also mention \cite{zang2020weak} as an alternative approach exploiting weak PDEs formulations, as well as \cite{li2020fourier}, which approximates operators by neural networks (however relying on training data from reference solutions), where a typical application is to map an initial condition to a solution of a PDE. For further references on approximating PDE solutions with neural networks we refer to the recent review articles \cite{beck2020overview,blechschmidt2021three,e2020algorithms,karniadakis2021physics}.
\par\bigskip

\textbf{Outline of the article.} In Section \ref{sec: PINN and BSDE loss} we review PINN and BSDE based approaches towards solving high-dimensional (parabolic) PDEs. In Section \ref{sec: diffusion loss} we introduce the diffusion loss, show its validity for solving PDEs of the form \eqref{eq: general boundary value problem}, and prove that it interpolates between the PINN and BSDE losses. Section \ref{sec: extensions} develops extensions of the proposed methodology to elliptic PDEs and eigenvalue problems. In Section \ref{sec: from losses to algorithms}, we discuss implementation details as well as some further modifications of the losses under consideration. In Section \ref{sec: numerics} we present numerical experiments, including a committor function example from molecular dynamics and a nonlinear eigenvalue problem motivated by quantum physics. Finally, Section \ref{sec: conclusion} concludes the paper with a summary and outlook.

\section{Variational formulations of boundary value problems}
\label{sec: PINN and BSDE loss}

In this section we consider boundary value problems such as \eqref{eq: general boundary value problem} in a variational formulation. That is, we aim at approximating the solution $V$ with some function $\varphi \in \mathcal{F}$ by minimizing suitable \emph{loss functionals} 
\begin{equation}
    \mathcal{L} : \mathcal{F} \to \R_{\ge 0},
\end{equation}
which are zero if and only if the boundary value problem is fulfilled,
\begin{equation}
    \mathcal{L}(\varphi) = 0 \iff \varphi = V.
\end{equation}
Here $\mathcal{F} \subset C^{2,1}(\Omega \times [0, T], \R) \cap C(\overline{\Omega} \times [0, T], \R)$
refers to an appropriate function class, usually consisting of deep neural networks. With a loss function at hand we can apply gradient-descent type algorithms to minimize (estimator versions of) $\mathcal{L}$, keeping in mind that different choices of losses lead to different statistical and computational properties and therefore potentially to different convergence speeds and robustness behaviors \cite{nusken2021solving}. 
\\

Throughout, we will work under the following assumption:
\begin{assumption}
\label{as: PDE}
The following hold:
\begin{enumerate}
\item 
The domain $\Omega$ is either bounded with piecewise smooth boundary, or $\Omega = \mathbb{R}^d$. 
\item 
The boundary value problem \eqref{eq: general boundary value problem} admits a unique classical solution $V \in C^{2, 1}(\Omega \times [0, T], \R) \cap C(\overline{\Omega} \times [0, T], \R)$. Moreover, the gradient of $V$ satisfies a polynomial growth condition in $x$, that is,
\begin{equation}
    |\nabla V(x,t)| \le C (1 + |x|^q), \qquad \qquad (x,t) \in \Omega \times[0,T],
\end{equation}
for some 
$C, q > 0$.
\item Complemented with a deterministic initial condition, the SDE \eqref{eq: uncontrolled SDE} admits a unique strong solution, globally in time.
\end{enumerate}
\end{assumption}

We note that the second assumption is trivially satisfied in the case when $\Omega$ is bounded.

\subsection{The PINN loss}

Losses based on PDE residuals go back to \cite{lagaris1998artificial, lagaris2000neural} and have gained recent popularity under the name \textit{physics-informed neural networks} (PINNs, \cite{raissi2017physics})  or \textit{deep Galerkin methods} (DGMs, \cite{sirignano2018dgm}). The idea is to minimize an appropriate $L^2$-error between the left- and right-hand sides of \eqref{eq: general PDE}-\eqref{eq: space boundary condition}, replacing $V$ by its approximation $\varphi$. The derivatives of $\varphi$ are computed analytically or via automatic differentiation and the data on which $\varphi$ is evaluated is distributed according to some prescribed probability measure (often a uniform distribution). A precise definition is as follows:
 
\begin{definition}[PINN loss]
\label{def:PINN}
Let $\varphi \in \mathcal{F}$. The \textit{PINN loss} consists of three terms,
\begin{equation}
\label{eq: PINN loss}
 \mathcal{L}_\text{PINN}(\varphi) = \alpha_{\text{int}} \mathcal{L}_\text{PINN,int}(\varphi) + \alpha_{\text{T}} \mathcal{L}_\text{PINN,T}(\varphi) + \alpha_{\text{b}} \mathcal{L}_\text{PINN,b}(\varphi),
 \end{equation}
 where
 \begin{subequations}
 \begin{align}
 \label{eq:PINN int}
     \mathcal{L}_\text{PINN,int}(\varphi) &= \E\left[\left( (\partial_t + L) \varphi(X, t) + h(X, t, \varphi(X, t), \sigma^\top\nabla \varphi(X, t)) \right)^2 \right], \\
     \label{eq:PINN T}
    \mathcal{L}_{\text{PINN,T}}(\varphi) &= \E\left[\left( \varphi(X^{(T)}, T) - f(X^{(T)}) \right)^2\right], \\
    \label{eq:PINN b}
    \mathcal{L}_{\text{PINN,b}}(\varphi) &= \E\left[\left( \varphi(X^{\text{b}}, t^\text{b}) - g(X^\text{b}, t^\text{b}) \right)^2\right].
\end{align}
\end{subequations}
Here, $\alpha_\text{int}, \alpha_\text{T}, \alpha_\text{b} > 0$ are suitable weights balancing the \emph{interior}, \emph{terminal} and \emph{boundary} constraints, and $(X,t) \sim \nu_{\Omega \times [0,T]}$, $X^{(T)} \sim \nu_{\Omega}$ and $(X^{\text{b}},t^{\text{b}}) \sim \nu_{\partial \Omega \times [0,T]}$
are distributed according to  probability measures $\nu_{\Omega \times [0,T]} \in \mathcal{P}(\Omega \times [0,T])$,  $\nu_{\Omega} \in \mathcal{P}(\Omega)$ and $\nu_{\partial \Omega \times [0,T]} \in \mathcal{P}(\partial \Omega \times [0,T])$ that are fully supported on their respective domains. 
\end{definition}

While  uniform distributions are canonical choices for $\nu_{\Omega \times [0,T]}$, $\nu_\Omega$ and $\nu_{\partial \Omega \times [0,T]}$, further research might reveal promising (possibly adaptive) alternatives that focus the sampling on specific regions of interest.

\begin{remark}[Generalizations]
Clearly, the loss contributions \eqref{eq:PINN int}-\eqref{eq:PINN b} represent in one-to-one correspondence the constraints in \eqref{eq: general PDE}-\eqref{eq: space boundary condition}.
By that principle, the PINN loss 
can straightforwardly be generalized to other types of PDEs, see \cite{karniadakis2021physics} and references therein.
Let us further already mention that choosing appropriate weights $\alpha_\text{int}, \alpha_\text{T}, \alpha_\text{b} > 0$ is crucial for algorithmic performance, but not straightforward. We will elaborate on this aspect in \Cref{sec: from losses to algorithms}.
\end{remark}

\begin{remark}[Unbounded domains] In the case when $\Omega = \mathbb{R}^d$, the boundary contribution \eqref{eq:PINN b} becomes obsolete and we set $\alpha_{\text{b}} = 0$. Note that in this scenario, it is impossible to choose $\nu_{\Omega \times [0,T]}$ and $\nu_{\Omega}$ to be uniform, and a choice to prioritize certain regions of $\Omega$ needs to be made. Although in theory the loss function \eqref{eq: PINN loss} still admits the solution to \eqref{eq: general PDE} as its unique minimizer, the lack of uniform coverage will make it challenging to obtain accurate approximations on sets which have small probability under $\nu_{\Omega \times [0,T]}$ and $\nu_{\Omega}$. Solving PDEs on large (and, in principle, unbounded) domains using the PINN loss hence requires a commensurate computational budget to ensure sufficient sampling for the expectations in \eqref{eq:PINN int}-\eqref{eq:PINN b}.
Analogous remarks apply to the BSDE and diffusion losses introduced below.
\end{remark}

\subsection{The BSDE loss}
\label{sec:BSDE}

The BSDE loss makes use of a stochastic representation of the boundary value problem \eqref{eq: general boundary value problem} given by a backward stochastic differential equation (BSDE) that is rooted in the correspondence between the differential operator $L$ defined in \eqref{eq: infinitesimal generator variational} and the stochastic process $X$ defined in \eqref{eq: uncontrolled SDE}. Indeed, according to \cite{pardoux1998backward}, the PDE \eqref{eq: general boundary value problem} is related to the system
\begin{subequations}
\label{eq: FBSDE system}
\begin{align}
\label{eq:forward SDE}
\mathrm{d} X_s &= b(X_s, s) \, \mathrm{d} s + \sigma(X_s, s) \, \mathrm{d} W_s,  \qquad   &X_{t_0} = x_{\mathrm{init}},  \\
\label{eq: BSDE}
\mathrm{d} Y_s &= -h(X_s, s, Y_s, Z_s) \,\mathrm{d}s + Z_s \cdot \mathrm{d} W_s, \qquad  &Y_{T \wedge \tau}  = k(X_{T\wedge\tau}, T\wedge\tau),
\end{align}
\end{subequations}
where $\tau = \inf \{t> 0 : X_t \notin \Omega \}$ is the first exit time from $\Omega$ and $k$ subsumes the boundary  conditions,
\begin{equation}
\label{eq: def k BSDE}
k(x,t) = \begin{cases} 
f(x), & t = T, \,\, x \in \Omega, \\
g(x, t), & t \le T,\,\, x \in \partial \Omega.
\end{cases}
\end{equation}
Owing to the fact that $X$ in \eqref{eq:forward SDE} is constrained at initial time $t_0$ and $Y$ in \eqref{eq: BSDE} at final time $T$, the processes $(X_s)_{t_0 \le s \le T}$ and $(Y_s)_{t_0 \le s \le T}$ are referred to as forward and backward, respectively.  
Given appropriate growth and regularity conditions on the coefficients $b$, $\sigma$, $h$ and $k$, It{\^o}'s formula implies that the backward processes satisfy
\begin{equation}
\label{eq: def Y, Z}
    Y_s = V(X_s, s), \qquad Z_s = \sigma^\top \nabla V(X_s, s),
\end{equation}
that is, they provide the solution to \eqref{eq: general boundary value problem} and its derivative along the trajectories of the forward process $(X_s)_{t_0 \le s \le T}$, see \cite{zhang2017backward}. Aiming for the approximation $\varphi \approx V$, we substitute \eqref{eq: def Y, Z} into \eqref{eq: BSDE} and replace $V$ by $\varphi$  to obtain
\begin{equation}
\widetilde{Y}_{T \wedge \tau} (\varphi) = \varphi(X_{t_0},t_0)  - \int_0^{T \wedge \tau} h(X_s,s,\varphi(X_s,s),\sigma^\top \nabla \varphi(X_s,s)) \, \mathrm{d}s + \int_0^{T \wedge \tau} \sigma^\top \nabla \varphi(X_s,s) \cdot \mathrm{d}W_s.
\end{equation}
Our notation $\widetilde{Y}_{T \wedge \tau} (\varphi)$ emphasizes the distinction from $Y_{T \wedge \tau}$ (which refers to the solution of \eqref{eq: FBSDE system}) and highlights the dependence on the particular choice $\varphi \in \mathcal{F}$. 
The key idea is now to penalize deviations from the terminal condition in \eqref{eq: BSDE} via the loss
\begin{equation}
    \mathcal{L}(\varphi) = \E \left[\left( k(X_{T\wedge\tau}, T\wedge\tau) - \widetilde{Y}_{T \wedge \tau}(\varphi) \right)^2 \right],
\end{equation}
see \cite{han2018solving}. We summarize this construction as follows.

\begin{definition}[BSDE loss]
\label{def: BSDE loss}
Let $\varphi \in \mathcal{F}$. The \textit{BSDE loss} is defined as
\begin{align}
\begin{split}
\label{eq: BSDE loss}
\mathcal{L}_\text{BSDE} (\varphi) = \mathbb{E} \Bigg[ &\Bigg(f(X_{\tau \wedge T})\mathbbm{1}_{\{\tau \wedge T = T\}} + g(X_{\tau \wedge T}, \tau \wedge T) \mathbbm{1}_{\{\tau \wedge T = \tau\}}  - \varphi(X_{t_0}, t_0)- \int_{t_0}^{\tau \wedge T} \sigma^\top \nabla\varphi(X_s, s) \cdot \mathrm{d}W_s \\
& \quad  + \int_{t_0}^{\tau \wedge T} h(X_s, s, \varphi(X_s, s), \sigma^\top \nabla\varphi(X_s, s)) \,  \mathrm ds\Bigg)^2 \Bigg],
\end{split}
\end{align}
where $(X_t)_{t_0 \le t \le \tau \wedge T}$ is a solution to \eqref{eq: uncontrolled SDE} and $\tau = \inf\{t > 0 : X_t \notin \Omega\}$ is the first exit time from $\Omega$. Furthermore, the initial condition $(X_0,t_0)$ is distributed according to a prescribed probability measure $\nu_{\Omega \times [0,T]}$ with full support on $\Omega \times [0,T]$.
\end{definition}

\begin{remark}[BSDE versus PINN]
\label{rem:comp BSDE PINN}
In contrast to the PINN loss \eqref{eq: PINN loss}, the BSDE loss \eqref{eq: BSDE loss} does not rely on a judicious tuning of the weights $\alpha_\text{int}, \alpha_\text{T}, \alpha_\text{b}$. On the other hand, implementations based on the BSDE loss face the challenge of simulating the hitting times $\tau = \inf\{t > 0 : X_t \notin \Omega\}$ efficiently and accurately. We shall elaborate on this aspect in \Cref{sec: SDE simulations}.  
\end{remark}

\begin{remark}[Neumann and periodic boundary conditions]
	Instead of the Dirichlet boundary condition \eqref{eq: space boundary condition}, Neumann or periodic boundary conditions may be considered, and the generalization of the PINN loss (as well as the diffusion loss introduced below) is straightforward. In the case of periodic boundary conditions, for instance, \eqref{eq:PINN b} may be replaced by
	\begin{equation}
	\label{eq: periodic boundary loss term}
	\mathcal{L}_{\text{b}}(\varphi) = \E\left[\left( \varphi(X^b) - \varphi(\overline{X}^b) \right)^2\right] + \E\left[\left| \nabla \varphi(X^b) - \nabla \varphi(\overline{X}^b) \right|^2\right],
	\end{equation}
	where $X^b \sim \nu_{\partial \Omega}$, and $\overline{X}^b$ refers to the reflected/periodic counterpart. 
	For the BSDE loss one could for instance incorporate homogeneous Neumann conditions by reflecting the stochastic process at the boundary \cite{han2020derivative}, however, straightforward extensions to more general boundary conditions seem to be not available. We refer to \cite{pardoux1998generalized}, which might provide a starting point for the construction of further related BSDE based methods.
\end{remark}

\begin{remark}[Relationship to earlier work]
\label{rem: BSDE loss related work}
The idea of approximating solutions to PDEs by solving BSDEs has been studied extensively \cite{bouchard2004discrete, gobet2005regression,pardoux1998backward}, where first approaches were regression based, relying on iterations backwards in time. These ideas do not seem to be straightforwardly applicable to the case when $\Omega$ is bounded, as the  trajectories of the forward process \eqref{eq:forward SDE} are not all of the same length (but see \cite{bouchard2009strong} and \cite{hartmann2019variational}).  A global variational strategy using neural networks has first been introduced in \cite{weinan2017deep}, where however in contrast to \Cref{def: BSDE loss}, the initial condition $(X_0,t_0)$ is deterministic (and fixed) and only parabolic problems on $\Omega= \mathbb{R}^d$ are considered. Moreover, slightly different choices for the approximations are chosen, namely $V$ is only approximated at $t_0=0$ and $\nabla V$ instead of $V$ is learnt by one neural network per time point (instead of using only one neural network with $t$ as an additional input).  
\end{remark}

\section{The diffusion loss}
\label{sec: diffusion loss}

In this section we introduce a novel loss that interpolates between the PINN and BSDE losses from Section \ref{sec: PINN and BSDE loss} using an auxiliary time parameter $\mathfrak{t} \in (0,\infty)$. As for the BSDE loss, the connection between the SDE \eqref{eq: uncontrolled SDE} and its infinitesimal generator \eqref{eq: infinitesimal generator variational} plays a major role:  It\^{o}'s formula
\begin{equation}
        V(X_T, T) - V(X_0, 0) = \int_0^T \left(\partial_s + L  \right)V(X_s, s) \,\mathrm ds + \int_0^T \sigma^\top \nabla V(X_s, s) \cdot  \mathrm dW_s
\end{equation}
motivates the following variational formulation of the  boundary value problem \eqref{eq: general boundary value problem}.

\begin{definition}[Diffusion loss]
\label{def: diffusion loss}
Let $\varphi \in \mathcal{F}$ and $\mathfrak{t} \in (0,\infty)$. The \textit{diffusion loss} consists of three terms,
\begin{equation}
\label{eq: definition diffusion loss}
    \mathcal{L}_{\text{diffusion}}^{\mathfrak{t}}(\varphi) = \alpha_{\text{int}} \mathcal{L}_{\text{diffusion,int}}^{\mathfrak{t}}(\varphi) + \alpha_{\text{T}} \mathcal{L}_{\text{diffusion,T}}^{\mathfrak{t}}(\varphi) + \alpha_{\text{b}} \mathcal{L}_{\text{diffusion,b}}^{\mathfrak{t}}(\varphi),
\end{equation}
where
\begin{subequations}
\begin{align}
\label{eq: definition diffusion loss int}
    \mathcal{L}_{\text{diffusion,int}}^{\mathfrak{t}}(\varphi) &= \E\Bigg[\Bigg(\varphi(X_{\mathcal{T}},  \mathcal{T}) - \varphi(X_{t_0}, t_0) - \int_{t_0}^{ \mathcal{T}}\sigma^\top\nabla \varphi(X_s, s) \cdot \mathrm dW_s  \\
    &\qquad\qquad + \int_{t_0}^{\mathcal{T}} h(X_s, s, \varphi(X_s, s), \sigma^\top\nabla \varphi(X_s, s))\,\mathrm ds\Bigg)^2\Bigg],\nonumber \\
    \label{eq: definition diffusion loss time boundary}
    \mathcal{L}_{\text{diffusion,T}}^{\mathfrak{t}}(\varphi) &= \E\left[\left( \varphi(X^{(T)}, T) - f(X^{(T)}) \right)^2\right], \\
    \label{eq: definition diffusion loss space boundary}
    \mathcal{L}_{\text{diffusion,b}}^{\mathfrak{t}}(\varphi) &= \E\left[\left( \varphi(X^{\text{b}}, t^{\text{b}}) - g(X^{\text{b}}, t^{\text{b}}) \right)^2\right],
\end{align}
 \end{subequations}
encode the constraints \eqref{eq: general PDE}-\eqref{eq: space boundary condition}, balanced by the weights
$\alpha_{\text{int}}, \alpha_{\text{T}}, \alpha_{\text{b}} > 0$. The process $(X_t)_{t_0 \le t \le \mathcal{T}}$ is a solution to \eqref{eq: uncontrolled SDE} with initial condition $(X_0,t_0)\sim \nu_{\Omega \times [0,T]}$ and maximal trajectory length $\mathfrak{t} > 0$. The stopping time  $\mathcal{T} := (t_0 + \mathfrak{t})\wedge \tau \wedge T$ is a shorthand notation, referring to the (random) final time associated to a realization of the path $X$ as it either hits the parabolic boundary $\partial \Omega \times \{T\}$ or reaches the maximal time $t_0 + \mathfrak{t}$. As in Definition \ref{def:PINN}, $\tau = \inf\{t > 0 : X_t \notin \Omega\}$ is the exit time from $\Omega$, and 
$(X^{\text{b}},t^{\text{b}}) \sim \nu_{\partial \Omega \times [0,T]}$, $X^{(T)} \sim \nu_{\Omega}$
are distributed according to probability measures that are fully supported on their respective domains.
\end{definition}

\begin{figure}[H]
\centering
\includegraphics[width=1.0\linewidth]{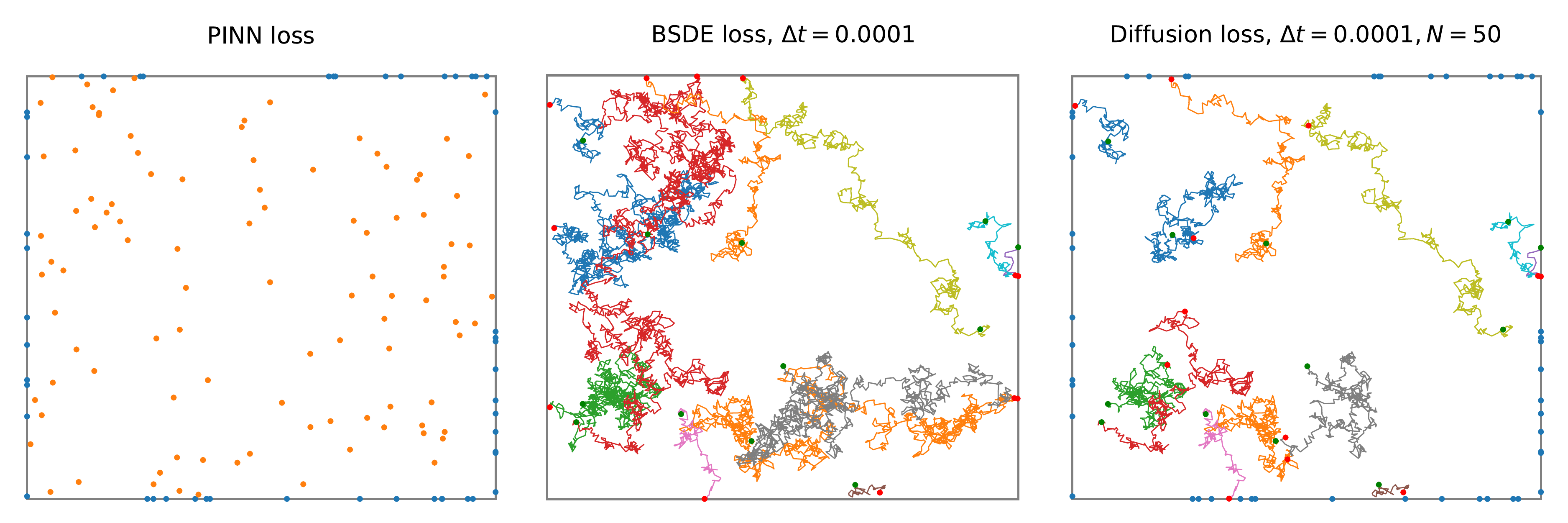}
\caption{We illustrate the training data used for the three losses inside the unit square $\Omega = (0,1)^2$. The PINN loss in the left panel  takes i.i.d. data points that are sampled from prescribed probability distributions in the domain $\Omega$ and on the boundary $\partial \Omega$ (in this case from uniform distributions). The BSDE loss (middle   panel) uses trajectories associated to the SDE \eqref{eq: uncontrolled SDE} that are started at random points $X_0$ (green points) and run until they hit the boundary (red points). The trajectories for the diffusion loss have a maximal length $\mathfrak{t}$ and therefore frequently start and end inside $\Omega$, as displayed in the right panel. The blue points for the PINN and diffusion losses indicate the additionally sampled boundary data.}
\label{fig: illustration of three losses}
\end{figure}

\begin{remark}[Comparison to the BSDE and PINN losses]
\label{rem:comparison BSDE and PINN}
In contrast to the PINN loss from Definition \ref{def:PINN}, the data inside the domain $\Omega$ is not sampled according to a prescribed probability measure $\nu_{\Omega}$, but along trajectories of the diffusion \eqref{eq: uncontrolled SDE}. Consequently, second derivatives of $\varphi$ do not have to be computed explicitly, but are approximated using the driving Brownian motion (and -- implicitly -- It{\^o}'s formula). A main difference to the BSDE loss from Definition \ref{def: BSDE loss} is that the simulated trajectories have a maximal length $\mathfrak{t}$, which might be beneficial computationally if the final time $T$ or the exit time $\tau$ is large (with high probability). Additionally, the sampling of extra boundary data circumvents the problem of accurately simulating those exit times (see Remark \ref{rem:comp BSDE PINN}). Both aspects will be further discussed in \Cref{sec: from losses to algorithms}. We refer to \Cref{fig: illustration of three losses} for a graphical illustration of the data required to compute the three losses.
\end{remark}

\begin{remark}[Reinforcement learning]
\label{rem:reinforcement}
The form of \eqref{eq: definition diffusion loss int} can directly be related to (least squares) temporal difference objectives \cite[Section 3.2]{szepesvari2010algorithms}, the PDE \eqref{eq: general PDE} being a continuous-time analog to the Bellman equation that governs optimal control. We refer to \cite[Section 2]{zhou2021actor} for a lucid exposition of the relevant ideas in the context of (elliptic) PDEs. Interestingly, however, the background on control and reinforcement learning (although providing further intuition) is not strictly necessary, as the diffusion loss applies to the 
 general parabolic PDE in \eqref{eq: general PDE}.\footnote{This PDE contains control-related Hamilton-Jacobi-Bellman PDEs as special cases, see, for instance \cite{pham2009continuous}.} The diffusion loss can hence be thought of as a generalization of the approach in \cite{zhou2021actor}, see in particular equation (2.27) therein, with $\gamma = 0$.
\end{remark}

The following proposition shows that the loss $\mathcal{L}_{\text{diffusion}}^{\mathfrak{t}}$ is indeed suitable for the boundary value problem \eqref{eq: general boundary value problem}.

\begin{proposition}
\label{prop: PDE fulfilled iff diffusion loss zero}
Consider the diffusion loss as defined in \eqref{eq: definition diffusion loss}, and assume that $b$ and $\sigma$ are globally Lipschitz continuous in $x$, uniformly in $t \in [0,T]$. Furthermore, assume the  following Lipschitz and boundedness conditions on $f$, $g$ and $h$,
\begin{align*}
     |f(x)| & \le C(1+ |x|^p), \\
     |g(x, t)| & \le C(1+ |x|^p), \\
     |h(t,x,y,z)| & \le C (1 + |x|^p + |y| + |z|),
     \\
    |h(t,x,y,z) - h(t,x,y',z) | & \le C | y - y'|,
     \\
    | h(t,x,y,z) - h(t,x,y,z')| &\le C |z - z'|,
\end{align*}
for appropriate constants $C, p \ge 0$ and all $x,y,z \in \Omega$, $ t 
\in [0,T]$. Finally, assume that Assumption \ref{as: PDE} is satisfied.
Then for $\varphi \in \mathcal{F}$ the following are equivalent:
\begin{enumerate}
    \item The diffusion loss vanishes on $\varphi$,
\begin{equation}
    \mathcal{L}_{\mathrm{diffusion}}^{\mathfrak{t}}(\varphi) = 0.
\end{equation}
\item $\varphi$ fulfills the boundary value problem \eqref{eq: general boundary value problem}.
\end{enumerate}
\end{proposition}

\begin{proof}
Denoting by $X_s$ the unique strong solution to \eqref{eq: uncontrolled SDE}, an application of It\^{o}'s lemma to $\varphi(X_s, s)$ yields
\begin{equation}
    \varphi(X_{\mathcal{T}}, \mathcal{T}) = \varphi(X_{t_0}, t_0) + \int_{t_0}^{\mathcal{T}} (\partial_s + L)\varphi(X_s, s)\,\mathrm ds + \int_{t_0}^{\mathcal{T}} \sigma^\top\nabla \varphi(X_s, s) \cdot \mathrm dW_s,
\end{equation}
almost surely.
Assuming that $\varphi$ fulfills the PDE \eqref{eq: general PDE}, it follows from the definition in \eqref{eq: definition diffusion loss int} that $\mathcal{L}_{\text{diffusion,int}}^{\mathfrak{t}}(\varphi) = 0$. Similarly, the boundary conditions \eqref{eq: time boundary condition} and \eqref{eq: space boundary condition} imply that $ \mathcal{L}_{\text{diffusion,T}}^{\mathfrak{t}}(\varphi) =  \mathcal{L}_{\text{diffusion,b}}^{\mathfrak{t}}(\varphi) = 0$. Consequently, we see that  $\mathcal{L}_{\text{diffusion}}^{\mathfrak{t}}(\varphi) = 0$.

For the converse direction, observe that $\mathcal{L}_{\text{diffusion}}^{\mathfrak{t}}(\varphi) = 0$ implies that 
\begin{equation}
\label{eq: diffusion loss expression}
    \varphi(X_{\mathcal{T}},  \mathcal{T}) = \varphi(X_{t_0}, t_0) + \int_{t_0}^{ \mathcal{T}}\sigma^\top\nabla \varphi(X_s, s) \cdot \mathrm dW_s - \int_{t_0}^{\mathcal{T}} h(X_s, s, \varphi(X_s, s), \sigma^\top\nabla \varphi(X_s, s))\, \mathrm ds, 
\end{equation}
almost surely, and that the same holds with $\varphi$ replaced by $V$. We proceed by defining the processes $\widetilde{Y}_s := \varphi(X_s,s)$ and  $\widetilde{Z}_s := \sigma^\top  \nabla \varphi (X_s,s)$, as well as $Y_s := V(X_s,s)$ and $Z_s := \sigma^\top  \nabla V (X_s,s)$. 
By the assumptions on $\varphi$, $b$ and $\sigma$, the processes $Y$, $Z$, $\widetilde{Y}$ and  $\widetilde{Z}$ are progressively measurable with respect to the filtration generated by $(W_t)_{t \ge 0}$ and moreover square-integrable. Furthermore, the relation \eqref{eq: diffusion loss expression} shows that the pairs $(Y,Z)$ and $(\widetilde{Y},\widetilde{Z})$ satisfy a BSDE with terminal condition $\xi := \varphi(X_\mathcal{T},\mathcal{T})$  on the random time interval $[t_0, \mathcal{T}]$. Well-posedness of the BSDE (see \cite[Theorems 1.2 and 3.2]{pardoux1998backward}) implies that $Y = \widetilde{Y}$ and $Z = \widetilde{Z}$, almost surely. Conditional on $t_0$ and $X_{t_0}$, we also have $V(X_{t_0},t_0) = Y^{X_{t_0},t_0} = \widetilde{Y}^{X_{t_0},t_0} = \varphi(X_{t_0},t_0)$, where the superscripts denote conditioning on the initial time $t_0$ and corresponding initial condition $X_{t_0}$, see \cite[Theorems 2.4 and 4.3]{pardoux1998backward}. Hence, we conclude that $\varphi = V$, $\nu_{\Omega} \times [0,T]$-almost surely, and the result follows from the continuity of $\varphi$ and $V$ and the assumption that $\nu_{\Omega \times [0,T]}$ has full support.  
\end{proof}

We have noted before that the diffusion loss combines aspects from the BSDE and PINN losses. In fact, it turns out that the diffusion loss can be interpreted as a specific type of interpolation between the two. The following proposition makes this observation precise.

\begin{proposition}[Relation of the diffusion loss to the PINN and BSDE losses]
\label{prop: Relation of diffusion loss to PINN and BSDE losses}
Let $\varphi \in \mathcal{F}$. Assuming that the measures $\nu_{\Omega \times [0,T]}$ in Definitions \ref{def:PINN} and \ref{def: diffusion loss} coincide, we have that 
\begin{equation}
\label{eq: diffusion loss to PINN loss}
    \frac{ \mathcal{L}_{\mathrm{diffusion,int}}^{\mathfrak{t}}(\varphi)}{\mathfrak{t}^2} \to \mathcal{L}_\mathrm{PINN,int}(\varphi),
\end{equation}
as $\mathfrak{t} \to 0$.
Moreover, if $\nu_{\Omega \times [0,T]}$ refers to the same measure in Definitions \ref{def:PINN} and \ref{def: BSDE loss}, then 
\begin{equation}
\label{eq: diffusion loss to BSDE loss}
   \mathcal{L}_{\mathrm{diffusion,int}}^{\mathfrak{t}}(\varphi) \to \mathcal{L}_\mathrm{BSDE} (\varphi),
\end{equation}
as $\mathfrak{t} \to \infty$.
\end{proposition}
\begin{proof}
It\^{o}'s formula shows that $\mathcal{L}_{\text{diffusion,int}}^{\mathfrak{t}}$ can be expressed as
\begin{align}
\mathcal{L}_{\text{diffusion,int}}^{\mathfrak{t}}(\varphi) &= \E\left[\left( \int_{t_0}^{\mathcal{T}} (\partial_s + L) \varphi(X_s, s)\,\mathrm ds + \int_{t_0}^{\mathcal{T}} h(X_s, s, \varphi(X_s, s), \sigma^\top\nabla \varphi(X_s, s))\,\mathrm ds \right)^2\right],
\end{align}
which implies the limit \eqref{eq: diffusion loss to PINN loss} by dominated convergence, noting that $\mathcal{T} \to t_0$ as $\mathfrak{t} \to 0$, almost surely. The relation \eqref{eq: diffusion loss to BSDE loss} follows immediately from the definition of $\mathcal{L}_\mathrm{BSDE}$ by noting that $\mathcal{T} \to \tau \wedge T$ as $\mathfrak{t} \to \infty$, almost surely.
\end{proof}

\section{Extensions to elliptic PDEs  and eigenvalue problems}
\label{sec: extensions}

In this section we show that the ideas reviewed and developed in the previous sections in the context of the parabolic boundary value problem \eqref{eq: general boundary value problem} can straightforwardly be extended to the treatment of elliptic PDEs and certain eigenvalue problems (and in fact, the diffusion loss for parabolic problems could have been derived using \cite[Section 2.2.3]{zhou2021actor} as a starting point). To begin with, it is helpful to notice that the boundary value problem \eqref{eq: general boundary value problem} can be written in the following slightly more abstract form:

\begin{remark}[Compact notation]
\label{rem: compact notation of PDE}
Consider the operator $\mathcal{A} := \partial_t + L$, the space-time domain $\Omega_{xt} := \Omega \times [0, T)$ with (forward) boundary\footnote{We notice in passing that the (ordinary) topological boundary of $\Omega_{xt}$ is given by $\partial \Omega_{xt} = \Omega \times \{0,T\} \cup \partial \Omega \times [0, T]$, so that $\partial \Omega_{xt} = \overrightarrow{\partial}\Omega_{xt} \cup \{0\} \times \Omega$.} $\overrightarrow{\partial} \Omega_{xt} := \Omega \times \{T\} \cup \partial \Omega \times [0, T]$, and the augmented variable $z = (x, t)^\top \in \Omega \times [0,T]$. Then problem \eqref{eq: general boundary value problem} can be presented as
\begin{subequations}
\label{eq:compact problem}
\begin{align}
    \mathcal{A} V(z) + h(z, V(z), \sigma^\top\nabla_x V(z)) &= 0,\qquad  && z \in \Omega_{xt}, \\
    V(z) &= k(z) ,\qquad  && z \in \overrightarrow{\partial} \Omega_{xt},
\end{align}
\end{subequations}
with $k$ defined as in \eqref{eq: def k BSDE}.
Relying on \eqref{eq:compact problem}, we can equivalently define the BSDE loss as
\begin{align}
\mathcal{L}_\text{BSDE} (\varphi) = \mathbb{E} \Bigg[ &\Bigg(k(X_{\tau_{xt}}, \tau_{xt})  - \varphi(X_{t_0}, t_0)- \int_{t_0}^{\tau_{xt}} \sigma^\top \nabla\varphi(X_s, s) \cdot \mathrm{d}W_s + \int_{t_0}^{\tau_{xt}} h(X_s, s, \varphi(X_s, s), \sigma^\top \nabla\varphi(X_s, s))  \,\mathrm ds\Bigg)^2 \Bigg],
\end{align}
where $\tau_{xt} = \inf\{t > 0 : X_t \notin \Omega_{xt}\}$ is the first exit time from $\Omega_{xt}$. The PINN and diffusion losses can similarly be rewritten in terms of the space-time domain $\Omega_{xt}$ and exit time $\tau_{xt}$.
\end{remark}

\subsection{Elliptic boundary value problems}
\label{sec: special PDEs}

Removing the time dependence from the solution (and from the coefficients $b$ and $\sigma$ determining $L$) we obtain the elliptic boundary value problem
\begin{subequations}
\label{eqn: elliptic PDE}
\begin{align}
\label{eqn: elliptic PDE - domain}
    L V(x) + h(x, V(x), \sigma^\top\nabla V(x)) &= 0,\qquad  && x \in \Omega , \\
\label{eqn: elliptic PDE - boundary}
    V(x) & = g(x), \qquad  && x \in \partial \Omega,
\end{align}
\end{subequations}
with the nonlinearity $h \in C(\R^d \times \R \times \R^{d}, \R)$. In analogy to \eqref{eq: FBSDE system}, the corresponding backward equation is given by
\begin{equation}
    \mathrm d Y_s = -h(X_s, Y_s, Z_s) \,\mathrm ds +  Z_s \cdot \mathrm dW_s, \qquad Y_\tau = g(X_\tau), 
\end{equation}
where $\tau = \{t > 0: X_t \notin \Omega \}$ is the first exit time from $\Omega$. Given suitable boundedness and regularity assumptions on $h$ and assuming that $\tau$ is almost surely finite, one can show existence and uniqueness of solutions $Y$ and $Z$, which, as before, represent the solution $V$ and its gradient along trajectories of the forward process \cite[Theorem 4.6]{pardoux1998backward}.

Therefore, the BSDE, PINN and diffusion losses can be applied to \eqref{eqn: elliptic PDE} with minor modifications: Owing to the fact that there is no terminal condition, we set $f=0$ in \eqref{eq: BSDE loss}, as well as $\alpha_{\text{T}} = 0$ in \eqref{eq: PINN loss} and \eqref{eq: definition diffusion loss}, making \eqref{eq:PINN T} and \eqref{eq: definition diffusion loss time boundary} obsolete. With the same reasoning, we set $T = \infty$, incurring $\tau \wedge T = \tau$ and $\mathcal{T} = (t_0 + \mathfrak{t})\wedge \tau$; these simplifications are relevant for the expressions \eqref{eq: BSDE loss} and \eqref{eq: definition diffusion loss int}.
Proposition \ref{prop: PDE fulfilled iff diffusion loss zero} and its proof can straightforwardly be generalized to the elliptic setting.
An algorithm for solving elliptic PDEs of the type \eqref{eqn: elliptic PDE} in the spirit of the BSDE loss has been suggested in \cite{kremsner2020deep}, using the same approximation framework as in \cite{weinan2017deep} (cf. \Cref{rem: BSDE loss related work}). We note that the solutions to linear elliptic PDEs often admit alternative variational characterizations in terms of energy functionals \cite{weinan2018deep}. An approach using the Feynman-Kac formula has been considered in \cite{grohs2020deep}.

\subsection{Elliptic eigenvalue problems}
\label{sec: eigenvalue problems}

We can extend the algorithmic approaches from Sections \ref{sec: PINN and BSDE loss} and \ref{sec: diffusion loss} to eigenvalue problems of the form
\begin{subequations}
\label{eq: eigenvalue problem}
\begin{align}
    L V(x) &= \lambda V(x),\qquad  && x \in \Omega, \\ 
    V(x) & = 0, \qquad  && x \in \partial \Omega,
\end{align}
\end{subequations}
corresponding to the choice $h(x, y, z) = - \lambda y$ in the elliptic PDE \eqref{eqn: elliptic PDE}. Note, however, that $h$ now depends on the unknown eigenvalue $\lambda \in \R$. Furthermore, we can consider nonlinear eigenvalue problems,
\begin{subequations}
\label{eq:nonlinear eigenvalue problem}
\begin{align}
    L V(x)+ h(x, V(x), \sigma^\top \nabla V(x)) &= \lambda V(x),\qquad  && x \in \Omega, \\ 
    V(x) & = 0, \qquad  && x \in \partial \Omega,
\end{align}
\end{subequations}

with a general nonlinearity $h \in C(\R^d \times \R \times \R^d, \R)$.\par\bigskip

For the linear problem \eqref{eq: eigenvalue problem} it is known that, given suitable boundedness and regularity assumption on $b$ and $\sigma$, there exists a unique principal eigenvalue with strictly positive eigenfunction in $\Omega$, see \cite[Theorem 2.3]{berestycki1994principal}. This motivates us to consider the above losses, now depending on $\lambda$, as well as enhanced with an additional term, preventing the trivial solution $V \equiv 0$.  We define
\begin{equation}
\label{eq: eigenfunction loss}
    \mathcal{L}^\mathrm{eigen}(\varphi, \lambda) = \mathcal{L}_\lambda(\varphi) +  \alpha_{\text{c}} \mathcal{L}_{\text{c}}(\varphi),
\end{equation}
where $ \mathcal{L}_\lambda(\varphi)$ stands for either the PINN, the BSDE, or the diffusion loss (with the  nonlinearity $h$ depending on $\lambda$), $\mathcal{L}_{\text{c}}(\varphi) = (\varphi(x_{\text{c}}) - 1)^2$, and $\alpha_{\text{c}} > 0$ is a weight. Here $x_{\text{c}} \in \Omega$ is chosen deterministically, preferably not too close to the boundary $\partial \Omega$. Clearly, the term $\mathcal{L}_{\text{c}}$ encourages $\varphi(x_{\text{c}}) = 1$, thus discouraging $\varphi \equiv 0$. We note that $V(x_\text{c}) = 1$ can be imposed on solutions to \eqref{eq: eigenvalue problem} without loss of generality, since the eigenfunctions are determined up to a multiplicative constant only. Avoiding $\phi \equiv 0$ for nonlinear eigenvalue problems of the form \eqref{eq:nonlinear eigenvalue problem} needs to be addressed on a case-by-case basis; we present an example in Section \ref{sec:nonlinear sch}. 
 \par\bigskip

The idea is now to minimize $\mathcal{L}^\mathrm{eigen}(\varphi, \lambda)$ with respect to  $\varphi \in \mathcal{F}$ and $\lambda \in \R$ simultaneously, while constraining the function $\varphi$ to be non-negative. According to the following proposition this is a valid strategy to determine the first eigenpair.

\begin{proposition}
\label{prop: eigenvalue loss}
Let $\Omega$ be bounded, and assume that $L$ is uniformly elliptic, that is, there exist constants $c_0,C_0 > 0$ such that 
\begin{equation}
c_0 | \xi |^2 \le \sum_{i,j=1}^d (\sigma \sigma^\top)(x) \xi_i \xi_j  \le C_0 |\xi|^2,
\end{equation}
for all $\xi \in \mathbb{R}^d$. Moreover, assume that $b$ is bounded.
Let $\varphi \in \mathcal{F}$ with $\varphi \ge 0$ and assume that $\mathcal{L}_{\lambda}(\varphi) = 0$ if and only if \eqref{eq: eigenvalue problem} is satisfied. Then the following are equivalent:
\begin{enumerate}
    \item 
    \label{it:eig1}
    $\varphi$ is the principal eigenfunction for \eqref{eq: eigenvalue problem} with principal eigenvalue $\lambda$ and normalization $\varphi(x_c) = 1$.
    \item 
    \label{it:eig2}
    $\mathcal{L}^\mathrm{eigen}(\varphi, \lambda)$ vanishes on 
    the pair $(\varphi,\lambda)$, 
    \begin{equation}
    \label{eq:eigen zero}
    \mathcal{L}^\mathrm{eigen}(\varphi, \lambda) = 0.
\end{equation}
\end{enumerate}
\end{proposition}

\begin{remark}
The assumption that $\mathcal{L}_\lambda(\varphi)$ is equivalent to \eqref{eq: eigenvalue problem} is satisfied for any `reasonable' loss function. For the diffusion loss, Proposition \ref{prop: PDE fulfilled iff diffusion loss zero} establishes this condition whenever the coefficients in \eqref{eq: eigenvalue problem} are regular enough.
\end{remark}

\begin{proof}
It is clear that \ref{it:eig1}. implies \ref{it:eig2}. by the  construction of \eqref{eq: eigenfunction loss}. For the converse direction, notice that \eqref{eq:eigen zero} implies $\varphi(x_c) = 1$ as well as \eqref{eq: eigenvalue problem}, that is, $\varphi$ is an eigenfunction with eigenvalue $\lambda$. In conjunction with the constraint $\varphi \ge 0$, it follows by \cite[Theorem 2.3]{berestycki1994principal} that $\varphi$ is the principal eigenfunction.
\end{proof}

An alternative approach towards \eqref{eq: eigenvalue problem} can be found in \cite{han2020solving}, where the eigenvalue problem is connected to a parabolic PDE and formulated as a fixed point problem.

\section{From losses to algorithms}
\label{sec: from losses to algorithms}

In this section we discuss some details regarding implementation aspects. 
For convenience, let us start by stating a prototypical algorithm  based on the losses introduced in \Cref{sec: PINN and BSDE loss,sec: diffusion loss}:\par\bigskip

\begin{algorithm}[H]
\SetAlgoLined
Choose a parametrization $\mathbb{R}^p \ni \theta \mapsto \varphi_{\theta}$.\\
Initialize $\varphi_\theta$ (with a parameter vector $\theta \in \R^p$). \\
Choose an optimization method $\texttt{descent}$, a batch size $K \in \mathbb{N}$ and a learning rate $\eta > 0$. For PINN and diffusion losses choose weights $\alpha_{\text{int}}, \alpha_{\text{b}}, \alpha_{\text{T}} > 0$ and batch sizes $K_{\text{b}}, K_{\text{T}} \in \mathbb{N}$. For BSDE and diffusion losses choose a step-size $\Delta t > 0$, for the diffusion loss choose a trajectory length $\mathfrak{t} > 0$.
\\
\Repeat{ convergence}{
Choose a loss function $\mathcal{L}$ from either \eqref{eq: PINN loss}, \eqref{eq: BSDE loss} or \eqref{eq: definition diffusion loss}.\\
Simulate data according to the chosen loss.\\
Compute $\widehat{\mathcal{L}}(\varphi_\theta)$ as a Monte Carlo version of $\mathcal{L}$.  \\
Compute $\nabla_{\theta}\widehat{\mathcal{L}}(\varphi_\theta)$ using automatic differentiation.\\
Update parameters: $\theta \gets \theta - \eta \, \texttt{descent}(\nabla_\theta\widehat{\mathcal{L}}(\varphi_\theta))$.
}
\caption{Approximation of the solution $V$ to the boundary value problem \eqref{eq: general boundary value problem}.}
\KwResult{$\varphi_\theta \approx V$.}
\label{algorithm}
\end{algorithm}\par\bigskip

\textbf{Function approximation.} In this paper, we rely on neural networks to provide the parametrization $\mathbb{R}^p \ni \theta \mapsto \varphi_{\theta}$ referred to in Algorithm \ref{algorithm} (but note that alternative function classes might offer specific benefits, see, for instance \cite{richter2021solving}). Standard feed-forward neural networks are given by
\begin{equation}
\label{eq: def NN}
\varphi_\theta:\R^d \to \R, \qquad \qquad
\varphi_\theta(x) = A_L \varrho(A_{L-1} \varrho(\cdots  \varrho(A_1 x + b_ 1) \cdots) + b_{L-1}) + b_L,
\end{equation}
with a collection of matrices $A_l \in \R^{n_{l} \times n_{l-1}}$ and vectors $b_l \in \R^{n_l}$ comprising the learnable parameters, $\theta = (A_l,b_l)_{1 \le l \le L}$. Here $L$ denotes the depth of the network, and we have $n_0 = d$ as well as $n_L  = 1$.
The nonlinear activation function $\varrho: \R \to \R$ is to be applied componentwise.

 Additionally, we define the \textit{DenseNet} \cite{huang2017, weinan2018deep} as a variant of the feed-forward neural network \eqref{eq: def NN}, containing additional skip connections,
\begin{equation}
\label{eq:densenet}
\varphi^{\mathrm{DenseNet}}_\theta(x) = A_L x_L + b_L,
\end{equation}
where $x_{L}$ is specified recursively as
\begin{align}
y_{l+1} = \varrho(A_l x_l + b_l), \qquad 
x_{l+1} = (x_l, y_{l+1})^\top, \qquad x_1 = x,
\end{align}
with $A_l \in \R^{n_l \times \sum_{i=0}^{l-1} n_i}$  and $b_l \in \R^l$ for $1 \le l \le L-1$, $n_0 = d$. Again, the collection of matrices $A_l$ and vectors $b_l$ comprises the learnable parameters.\par\bigskip

\textbf{Comparison of the losses (practical challenges).} The PINN, BSDE and diffusion losses differ in the way training data is generated (see Figure \ref{fig: illustration of three losses} and Table \ref{tab: loss comparisons data generation}); hence, the corresponding implementations face different challenges (see Table \ref{tab: loss comparisons challenges}). 

\begin{table}[H]
	\centering
	\captionsetup[subtable]{position=below}
	\captionsetup[table]{position=top}
	\caption{Comparison of the different losses.}
	\begin{subtable}[t]{0.44\linewidth}
		\centering
		\begin{tabular}{c|c|c|c}
			& \textbf{PINN} & \textbf{BSDE} & \textbf{Diffusion} \\ \hline
			SDE simulation & & \xmark & \xmark \\ \hline
			boundary data   & \xmark &  & \xmark\\ 
		\end{tabular}
		\caption{The three losses can be characterized by how training data is generated.}
		\label{tab: loss comparisons data generation}
	\end{subtable}%
	\hspace*{1em}
	\begin{subtable}[t]{0.56\linewidth}
		\centering
		\begin{tabular}{c|c|c|c}
			& \textbf{PINN} & \textbf{BSDE} & \textbf{Diffusion}  \\ \hline
			Hessian computations & \xmark &  &  \\ \hline
			exit time computations   &  & \xmark & \\ \hline
			weight tuning  & \xmark &  & \xmark\\ \hline
			long runtimes   &  & \xmark & \\  \hline
			discretization  &  & \xmark & \xmark \\ 
		\end{tabular}
		\caption{In this table we list potential challenges and drawbacks for the corresponding losses.}
		\label{tab: loss comparisons challenges}
	\end{subtable}
	\label{tab: loss comparisons}
\end{table}

First, the BSDE and diffusion losses rely on trajectorial data obtained from the SDE \eqref{eq: uncontrolled SDE}, in contrast to the PINN loss (cf. the first row in Table \ref{tab: loss comparisons data generation}). As a consequence, the BSDE and diffusion losses do not require the computation of second-order derivatives, as those are approximated implicitly using It{\^o}'s formula and the SDE \eqref{eq: uncontrolled SDE}, cf. the first row in Table \ref{tab: loss comparisons data generation}. From a computational perspective, the PINN loss therefore faces a significant overhead in high dimensions when the diffusion coefficient $\sigma$ is not sparse (as the expression \eqref{eq:PINN int} involves $d^2$ second-order partial derivatives). We notice in passing that an approach similar to the diffusion loss circumventing this problem has been proposed in 
\cite[Section 3]{sirignano2018dgm}.
On the other hand, evaluating the BSDE and diffusion losses requires discretizing the SDE \eqref{eq: uncontrolled SDE}, incurring additional numerical errors (cf. the last row in Table \ref{tab: loss comparisons challenges} and the discussion below in Section \ref{sec: SDE simulations}).  

Second, the PINN and diffusion losses incorporate boundary and final time constraints (see \eqref{eq: space boundary condition} and \eqref{eq: time boundary condition}) explicitly by sampling additional boundary data (see \eqref{eq:PINN T}, \eqref{eq:PINN b}, \eqref{eq: definition diffusion loss time boundary}, \eqref{eq: definition diffusion loss space boundary} and cf. the second row in Table \ref{tab: loss comparisons data generation}). On the one hand, this approach necessitates choosing the weights  $\alpha_{\text{int}}, \alpha_{\text{b}}, \alpha_{\text{T}} > 0$; it is by now well established that while algorithmic performance depends quite sensitively on a judicious tuning of these weights, general and principled guidelines to address this issue are not straightforward (see, however, \cite{van2020optimally, wang2021understanding, wang2020eigenvector, wang2022and}). Weight-tuning, on the other hand, is not required for implementations relying on the BSDE loss, as the boundary data is accounted for implicitly by the hitting event $\{(X_t,t) \notin \Omega \times [0,T) \}$ and the corresponding first two terms on the right-hand side of \eqref{eq: BSDE loss}. The hitting times $\tau = \inf \{t > 0 : X_t \notin \Omega\}$ may however be large, leading to a computational overhead in the generation of the training data (but see \ref{sec: forward control}), and are generally hard to compute accurately (but see Section \ref{sec: SDE simulations}).

\subsection{Simulation of diffusions and their exit times}
\label{sec: SDE simulations}

The BSDE and diffusion losses rely on trajectorial data obtained from the stochastic process defined in \eqref{eq: uncontrolled SDE}. In practice, we approximate this SDE on a time grid $t_0 \le t_1 \le \dots \le t_N$, for instance using the Euler-Maruyama scheme \cite{kloeden1992stochastic}

\begin{equation}
\label{eq: euler maruyama}
    \widetilde{X}_{n+1} = \widetilde{X}_n +  b(\widetilde{X}_n, t_n) \Delta t + \sigma(\widetilde{X}_n, t_n)\sqrt{\Delta t} \, \xi_{n+1},
\end{equation}
or, to be precise, by its stopped version
\begin{equation}
\label{eq: euler maruyama stopped}
    \widehat{X}_{n+1} = \widehat{X}_n +  \left( b(\widehat{X}_n, t_n) \Delta t + \sigma(\widehat{X}_n, t_n) \sqrt{\Delta t} \, \xi_{n+1} \right)\mathbbm{1}_{\mathcal{C}_{n+1}}
\end{equation}
with step condition $\mathcal{C}_{n} := \left\{\widetilde{X}_{n} \in \Omega\right\} \vee \left\{t_n \le T \right\}$ and time-increment $t_{n+1} = t_n + \Delta t \, \mathbbm{1}_{\mathcal{C}_{n+1}}$, where $\Delta t$ is the step-size and $\xi_{n+1}\sim\mathcal{N}(0, \mathrm{Id}_{d \times d})$ are standard normally distributed random variables. We can then straightforwardly construct Monte Carlo estimator versions of either the BSDE or the diffusion loss. For  example, the discrete version of the domain part of the diffusion loss \eqref{eq: definition diffusion loss int} reads
\begin{subequations}
\begin{align}
    \widehat{\mathcal{L}}_\mathrm{diffusion,int}^{(K, N)}(\varphi) &= \frac{1}{K}\sum_{k=1}^K \Bigg(\varphi(\widehat{X}_N^{(k)}, t_N^{(k)}) - \varphi(\widehat{X}_0^{(k)}, t_0^{(k)})  - \sum_{n=0}^{N-1}  \sigma^\top \nabla \varphi(\widehat{X}_n^{(k)}, t_n^{(k)}) \cdot \xi_{n+1}^{(k)} \sqrt{\Delta t} \, \mathbbm{1}_{\mathcal{C}_{n}^{(k)}} \\
    &\qquad\qquad\qquad + \sum_{n=0}^{N-1} h\left(\widehat{X}_n^{(k)}, t_n^{(k)}, \varphi(\widehat{X}_n^{(k)}, t_n^{(k)}), \sigma^\top \nabla\varphi(\widehat{X}_n^{(k)}, t_n^{(k)})\right)\Delta t \, \mathbbm{1}_{\mathcal{C}_{n}^{(k)}} \Bigg)^2,
\end{align}
\end{subequations}
where $K$ is the sample size, $N = \frac{\mathfrak{t}}{\Delta t}$ is the maximal discrete-time trajectory length, and $(\widehat{X}_n^{(k)})^{k=1...K}_{n=1...N}$ are independent copies of the iterates from \eqref{eq: euler maruyama stopped}. 
The Monte Carlo version of the BSDE loss can be formed analogously.
\par\bigskip

\begin{figure}[H]
\centering
\includegraphics[width=1.0\linewidth]{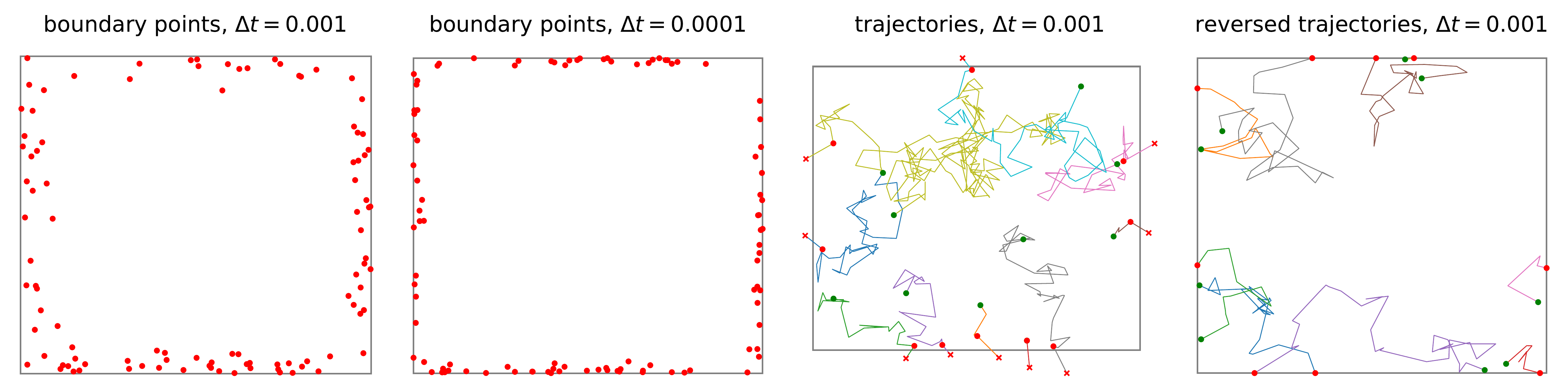}
\caption{Illustration of the boundary data in the BSDE method.}
\label{fig: BSDE boundary data illustration}
\end{figure}

Given suitable growth and regularity assumptions on the coefficients, the discretization errors of the forward and backward processes are of order $\sqrt{\Delta t}$ \cite{kloeden1992stochastic, zhang2017backward}, cf. also \cite{han2020convergence} for a numerical analysis on the original version of the BSDE loss. 
However, the stopped Euler-Maruyama scheme \eqref{eq: euler maruyama stopped} incurs additional errors due to the approximation of the first exit times from $\Omega$. In the two left panels of \Cref{fig: BSDE boundary data illustration} we illustrate this problem by displaying multiple ``last locations'' obtained according to \eqref{eq: euler maruyama stopped}, using two different step-sizes $\Delta t$. Clearly, all these points should in principle lie on the boundary, (naively) requiring a computationally costly small step-size.

Sophisticated approaches towards the accurate simulation of exit times for diffusions discretizing exit times have been put forward, see, e.g. \cite{buchmann2003solving, martin2022solving, hausenblas1999numerical}. For our purposes, however, it is not essential to compute the exit times, as long as the simulated trajectories are stopped accurately. We therefore suggest the following two attempts that aim at improving the sampling of boundary data:

\begin{enumerate}
    \item Rescaling: Start $X_0$ randomly in $\Omega$, simulate the trajectory and stop once the boundary has been crossed, however scale the last time step in such a way that the trajectory exactly ends on $\partial \Omega$.
    \item Time-reversal: Start $X_0$ on the boundary $\partial \Omega$ and simulate the trajectory for a given time $T$ (unless it hits the boundary again before time $T$, in this case stop the trajectory accordingly or resimulate). Then reverse the process such that the reversed process exactly ends on the boundary.
\end{enumerate}

An illustration of these two strategies can be found in the two right panels of \Cref{fig: BSDE boundary data illustration}. In our numerical experiments  we have found that both rescaling and time-reversal improves the performance of the BSDE method somewhat, but further research is needed.

\subsection{Further modifications of the losses}
\label{eq: further modifications}

In the following we discuss modifications of the PINN, BSDE and diffusion losses, relating also to versions that have appeared in the literature before.

\subsubsection{Forward control}
\label{sec: forward control}

We can modify the SDE-based BSDE and diffusion losses by including control terms $v \in C(\R^d \times [0, T], \R^d)$ in the forward process \eqref{eq: uncontrolled SDE}, yielding the controlled diffusion
\begin{equation}
\mathrm d X^v_s = \left(b(X^v_s,s) + \sigma(X^v_s,s) v(X^v_s,s)\right) \mathrm ds + \sigma(X^v_s,s) \, \mathrm dW_s.
\end{equation}
Applying It\^{o}'s formula we may obtain losses similar to those considered in Section \ref{sec: PINN and BSDE loss} and \ref{sec: diffusion loss}. For instance, alternative versions of the diffusion loss take the form
\begin{align}
\begin{split}
 \mathcal{L}_{\text{diffusion,int}}^{\mathfrak{t}, v}(\varphi) &= \E\Bigg[\Bigg(\varphi(X^v_{\mathcal{T}},  \mathcal{T}) - \varphi(X^v_{t_0}, t_0) - \int_{t_0}^{ \mathcal{T}}\sigma^\top\nabla \varphi(X_s, s) \cdot \mathrm dW_s  \\
    &\qquad\qquad + \int_{t_0}^{\mathcal{T}} \left[h(X_s^v, s, \varphi(X_s^v, s), \sigma^\top\nabla \varphi(X_s^v, s)) -  v(X^v_s,s) \cdot  \sigma^\top \nabla\varphi(X_s^v, s)\right]\mathrm ds\Bigg)^2\Bigg],
\end{split}
\end{align}
replacing \eqref{eq: definition diffusion loss int}.
We note in passing that Proposition \ref{prop: PDE fulfilled iff diffusion loss zero} extends straightforwardly to $\mathcal{L}_{\text{diffusion,int}}^{\mathfrak{t}, v}$ under the assumption that $v$ satisfies appropriate Lipschitz and growth conditions.
\par\bigskip

Similar considerations apply for the BSDE loss, noting that for solutions to the generalized BSDE system \cite{nusken2021solving}
\begin{subequations}
\begin{alignat}{2}
\mathrm d X^v_s &= \left(b(X^v_s,s) + \sigma(X^v_s,s) v(X^v_s,s)\right) \mathrm ds + \sigma(X^v_s,s) \, \mathrm dW_s, \qquad \qquad \qquad \qquad &X^v_{t_0} = x_\mathrm{init}, \\
\mathrm{d}Y_s^{v} &= -h(X^v_s, s, Y^v_s, Z^v_s) \, \mathrm{d}s + v(X^v_s,s) \cdot Z_s^v \, \mathrm{d}s + Z_s^v \cdot \mathrm{d}W_s, \qquad \qquad  &Y_T^{v} = k(X^v_{T\wedge\tau}, T\wedge\tau),
\end{alignat}
\end{subequations}
we still have the relations
\begin{equation}
    Y_s^v = V(X_s^v, s), \qquad Z_s^v = \sigma^\top \nabla V(X_s^v, s)
\end{equation}
for suitable $v \in C(\R^d \times [0, T], \R^d)$, analogously to \eqref{eq: def Y, Z}. This immediately incurs the family of losses
\begin{align}
\begin{split}
\mathcal{L}_\text{BSDE}^v (\varphi) = \mathbb{E} \Bigg[ &\Bigg(f(X^v_{\tau \wedge T})\mathbbm{1}_{\tau \wedge T = T} + g(X^v_{\tau \wedge T}, \tau \wedge T) \mathbbm{1}_{\tau \wedge T = \tau} - \varphi(X^v_{t_0}, t_0)- \int_{t_0}^{\tau \wedge T} \sigma^\top \nabla\varphi(X^v_s, s) \cdot \mathrm{d}W_s \\
&   \quad + \int_{t_0}^{\tau \wedge T} \left( h(X^v_s, s, \varphi(X_s, s), \sigma^\top \nabla\varphi(X_s, s))-  v(X^v_s,s) \cdot  \sigma^\top \nabla\varphi(X_s^v, s)  \right) \mathrm ds\Bigg)^2 \Bigg],
\end{split}
\end{align}
parametrized by $v \in C(\R^d \times [0, T], \R^d)$.\par\bigskip

Adding a control to the forward process can be understood as driving the data generating process into regions of interest, for instance possibly alleviating the problem that exit times might be large (see Table \ref{tab: loss comparisons challenges} and the corresponding discussion). Identifying suitable forward controls might be an interesting topic for future research (we refer to \cite{nusken2021solving} for some systematic approaches in this respect relating to Hamilton-Jacobi-Bellman PDEs).

\subsubsection{Approximating the gradient of the solution}

The constraints imposed by the BSDE system \eqref{eq: FBSDE system} can be enforced by losses that are slightly different from \Cref{def: BSDE loss}. Going back to \cite{weinan2017deep}, we can for instance use the fact that the backward process $Y$ can be written in a forward way, yielding the discrete-time process
\begin{equation}
\label{eq: discrete backward process}
    \widehat{Y}_{n+1} = \widehat{Y}_n - h(\widehat{X}_n, t_n, \widehat{Y}_n, \widehat{Z}_n) \Delta t + \widehat{Z}_n \cdot \xi_{n+1} \sqrt{\Delta t}.
\end{equation}
The scheme  \eqref{eq: discrete backward process} is explicit, the unknowns being $\widehat{Y}_0$ and $\widehat{Z}_n$, for $n \in \{0, \dots, N-1 \}$. This motivates approximating the single parameter $y_0 \approx \widehat{Y}_0 \in \R$ as well as the vector fields $ \phi \approx \sigma^\top \nabla V \in C(\R^d \times [0, T], \R^d)$, rather than $V$ directly. This approach gives rise to the loss
\begin{align}
\label{eq:Jentzen}
\begin{split}
\mathcal{L}_\mathrm{BSDE-2} (\phi, y_0) = \mathbb{E} \Bigg[& \Bigg(f(X_{\tau \wedge T})\mathbbm{1}_{\tau \wedge T = T} + g(X_{\tau \wedge T}, \tau \wedge T) \mathbbm{1}_{\tau \wedge T = \tau} - y_0 - \int_0^{\tau \wedge T} \phi(X_s, s) \cdot \mathrm{d}W_s \\
&\qquad + \int_0^{\tau \wedge T} h(X_s, s, Y_s, \phi(X_s, s))  \,\mathrm ds \Bigg)^2 \Bigg].
\end{split}
\end{align}
In this setting $X_0$ has to be chosen deterministically;  we note that a potential drawback is thus that the solution is only expected to be approximated accurately in regions that can be reached by the forward process $X_t$ (starting at $X_0$) with sufficiently high probability. It has been shown in \cite{nusken2021solving} that  alternative losses (like the log-variance loss) can be considered whenever the nonlinearity $h$ only depends on the solution through its gradient, in which case the extra parameter $y_0$ can be omitted.

\subsubsection{Penalizing deviations from the discrete scheme}

Another approach that is rooted in the discrete-time backward process \eqref{eq: discrete backward process} has been suggested in \cite{raissi2018forward} for problems on unbounded domains. It relies on the idea to penalize deviations from \eqref{eq: discrete backward process}, for each $n \in \{0, \dots, N-1 \}$ (cf. also \cite{hure2020deep, richter2021solving}, where however an implicit scheme and backward iterations are used). Aiming for $\widehat{Y}_n \approx \varphi(\widehat{X}_n, t_n)$, $\widehat{Z}_n \approx \sigma^\top \nabla \varphi(\widehat{X}_n, t_n)$, this motivates the loss
\begin{equation}
\label{eq: Raissi loss}
    \widehat{\mathcal{L}}^{(K, N)}_\mathrm{BSDE-3}(\varphi) = \alpha_{\text{int}}\widehat{\mathcal{L}}^{(K, N)}_\mathrm{BSDE-3,int}(\varphi) + \alpha_{\text{b}}\widehat{\mathcal{L}}^{(K, N)}_\mathrm{BSDE-3,b}(\varphi)
\end{equation}
with interior part
\begin{align}
\label{eq: Raissi loss interior}
    \widehat{\mathcal{L}}^{(K, N)}_\mathrm{BSDE-3,int}(\varphi) = \frac{1}{K} \sum_{k=1}^K\sum_{n=0}^{N-1}\left(\varphi(\widehat{X}_{n+1}^{(k)}) - \varphi(\widehat{X}_{n}^{(k)}) + h\left(\widehat{X}_{n}^{(k)},\varphi(\widehat{X}_{n}^{(k)}), \sigma^\top\nabla\varphi(\widehat{X}_{n}^{(k)}) \right) \Delta t - \sigma^\top \nabla\varphi(\widehat{X}_{n}^{(k)})\xi_{n+1} \sqrt{\Delta t} \right)^2
\end{align}
and boundary term
\begin{align}
    \widehat{\mathcal{L}}^{(K, N)}_\mathrm{BSDE-3,b}(\varphi) = \frac{1}{K} \sum_{k=1}^K \left(\varphi(\widehat{X}_{N}^{(k)}) - g(\widehat{X}_{N}^{(k)}) \right)^2.
\end{align}
A generalization to equations posed on bounded domains is straightforward. We also note that in contrast to the diffusion loss, \eqref{eq: Raissi loss} does not seem to naturally derive from a continuous-time formulation.
 Interestingly, $\widehat{\mathcal{L}}^{(K, N)}_\mathrm{BSDE-3}$  can be related to the diffusion loss via Jensen's inequality, 
\begin{equation}
    \widehat{\mathcal{L}}_\mathrm{diffusion,int}^{(K, N)}(\varphi) \le N \widehat{\mathcal{L}}^{(K, N)}_\mathrm{BSDE-3,int}(\varphi).
\end{equation}

Yet another approach that is based on a discrete backward scheme is the following. Let us initialize $\widehat{Y}_0 = \varphi(\widehat{X}_0, 0)$ and simulate
\begin{equation}
    \widehat{Y}_{n+1} = \widehat{Y}_n - h(\widehat{X}_n, \widehat{Y}_n, \sigma^\top \nabla \varphi(\widehat{X}_n, t_n))\Delta t + \sigma^\top \nabla \varphi(\widehat{X}_n, t_n) \cdot \xi_{n+1} \sqrt{\Delta t},
\end{equation}
for $n \in \{0, \dots, N-1 \}$, where, similarly to \eqref{eq:Jentzen}, but in in contrast to \eqref{eq: Raissi loss interior}, only $\widehat{Z}_n$ is replaced by its approximation $\sigma^\top \nabla \varphi(\widehat{X}_n, t_n)$ while $\widehat{Y}_n$ is retained from previous iteration steps. Again penalizing deviations from the discrete-time scheme, we can now introduce the loss
\begin{equation}
\label{eq: Jentzen loss}
    \widehat{\mathcal{L}}^{(K, N)}_\mathrm{BSDE-4}(\varphi) = \frac{\alpha_1}{K} \sum_{k=1}^K\sum_{n=0}^N\left(\varphi(\widehat{X}_n^{(k)}, t_n) - \widehat{Y}_n^{(k)}\right)^2 + \frac{\alpha_2}{K} \sum_{k=1}^K\left( \varphi(X_b^{(k)}) - g(X_b^{(k)}) \right)^2.
\end{equation}

Both in $\mathcal{L}_\mathrm{BSDE-3}$ and $\mathcal{L}_\mathrm{BSDE-4}$ the deterministic initial condition $\widehat{X}_0$ at $t = 0$ can be replaced by random choices $(\widehat{X}_{t_0}, t_0) \sim \nu_{\Omega \times [0,T]}$, adjusting the sums in \eqref{eq: Raissi loss interior} and \eqref{eq: Jentzen loss} accordingly.

\section{Numerical experiments}
\label{sec: numerics}

In this section we provide several numerical examples of high-dimensional parabolic and elliptic PDEs that shall demonstrate the performances of \Cref{algorithm} using the different loss functions introduced previously. We focus on the PINN, BSDE and diffusion losses from \Cref{sec: PINN and BSDE loss,sec: diffusion loss} since their modified versions from \Cref{eq: further modifications} in general led to similar or worse performances. If not specified otherwise, the approximation of  $\varphi$ relies on a DenseNet architecture defined in \eqref{eq:densenet} with ReLU activation function and four hidden layers with $d+20$, $d$, $d$, $d$ units respectively (recall that $d$ specifies the dimension). The optimization is carried out using the Adam optimizer \cite{kingma2014adam} with standard parameters and learning rate $\eta = 0.001$. Throughout, we take $K_{\text{int}}=200$ samples inside the domain $\Omega$ and (for the PINN and diffusion losses) $K_{\text{b}}=50$ samples on the boundary, per gradient step. For the SDE discretization we  choose a step-size of $\Delta t = 0.001$.
The weight configurations for the PINN and diffusion losses are optimized manually; we then only report the results for the optimal settings (see the discussion on weight-tuning in Section \ref{sec: from losses to algorithms}). We refer to the code at \url{https://github.com/lorenzrichter/path-space-PDE-solver}.

\subsection{Nonlinear toy problems}

Let us start with a nonlinear toy problem for which analytical reference solutions are available. Throughout this subsection, the domain  of interest is taken to be the unit ball $\Omega = \{ x \in \R^d : |x| < 1 \}$.

\subsubsection{Elliptic problem with Dirichlet boundary data}
\label{sec: Elliptic problem with Dirichlet boundary data}
We first consider an elliptic boundary value problem of the form \eqref{eqn: elliptic PDE}. Let $\gamma \in \R$ and choose
\begin{subequations}
\label{eq:toy}
  \begin{equation}
            b(x, t) = \mathbf{0}, \qquad \sigma(x, t) = \sqrt{2}\,  \mathrm{Id}_{d \times d}, \qquad g(x) = e^\gamma, \\
  \end{equation}
    \begin{equation}
      h(x, y, z) = -2\gamma y (\gamma |x|^2 + d) + \sin\left(e^{2\gamma |x|^2} - y^2\right).
  \end{equation}
  \end{subequations}
It is straightforward to verify that 
\begin{equation}
\label{eq: exponential on ball solution}
V(x) = e^{\gamma |x|^2}
\end{equation}
is the unique solution to \eqref{eqn: elliptic PDE}. \par\bigskip

\begin{figure}[H]
\centering
  \centering
  \includegraphics[width=1.0\linewidth]{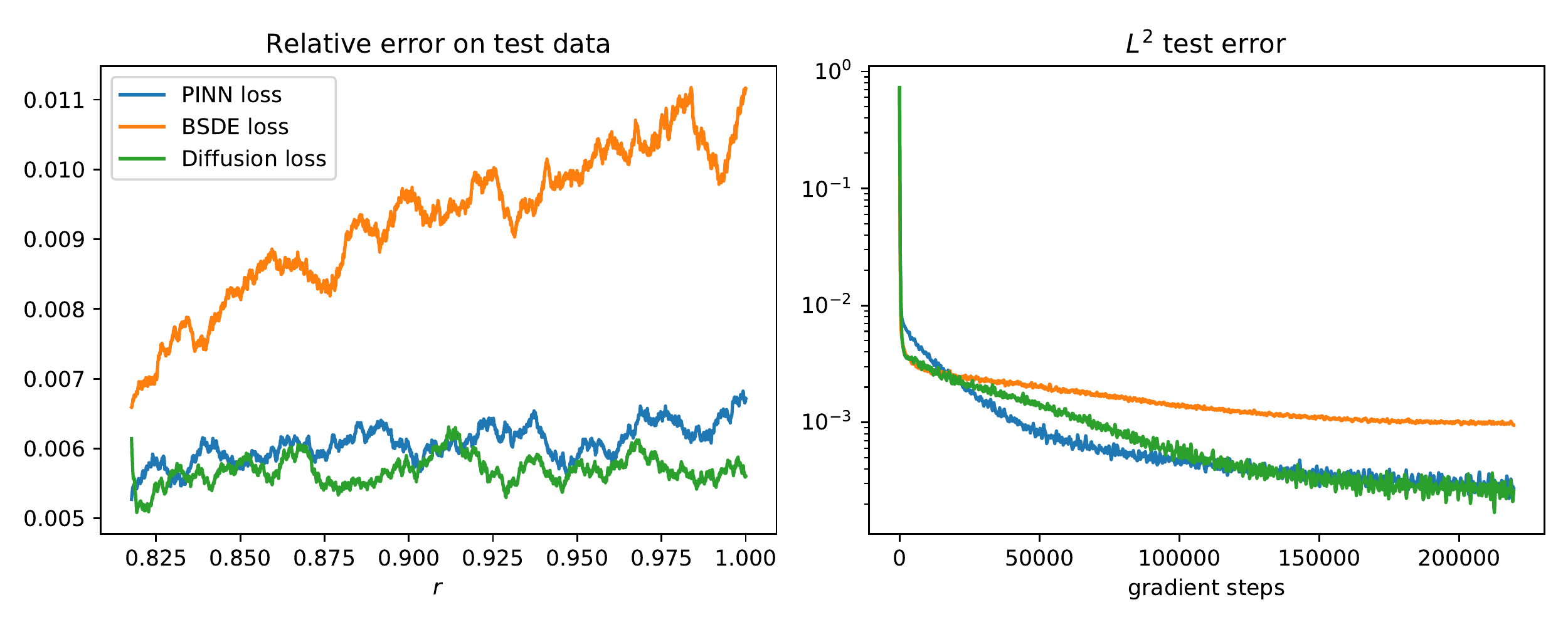}
  \caption{Left: Average relative errors as a function of $r = |x|$ evaluated on uniformly sampled data for the three losses smoothed with a moving average over $500$ data points. For the training the diffusion loss relies on trajectory lengths of $\mathfrak{t} = 0.02$. Right: $L^2$ error during the training iterations evaluated on uniformly sampled test data.}
  \label{fig: exponential on ball}
\end{figure}

We consider $d = 50$ and set $\gamma = 1$. For the PINN and diffusion losses the optimized weights are given by $\alpha_{\text{int}} = 10^{-5}$, $\alpha_{\text{b}} = 1$ and $\alpha_{\text{int}} = 0.1$, $\alpha_{\text{b}} = 1$, respectively. We sample the data uniformly and take a maximal (discrete-time) trajectory length of $N = 20$ (i.e. $\mathfrak{t} = 0.02$) for the diffusion loss. In \Cref{fig: exponential on ball} (left panel) we display the average relative errors $\frac{|\varphi(x) - V(x)|}{V(x)}$ as a function of $r = |x|$. Due to volume distortion in high dimensions, very few samples are drawn close to the center of the ball, and hence the numerical results appear to be unreliable for $r \le 0.8$. While this effect could be alleviated by changing the measure $\nu_{\Omega}$ accordingly, we content ourselves here with a comparison for $r \in [0.8,1]$.
In the right panel we display the $L^2$ error during the training iterations evaluated on uniformly sampled test data. We observe that the PINN and diffusion losses yield similar results and that the BSDE loss performs worse, in particular close to the boundary. We attribute this effect to the challenges inherent in the simulation of hitting times, see \Cref{sec: SDE simulations}.

\begin{figure}[H]
\centering
\begin{subfigure}{.5\textwidth}
  \centering
  \includegraphics[width=\linewidth]{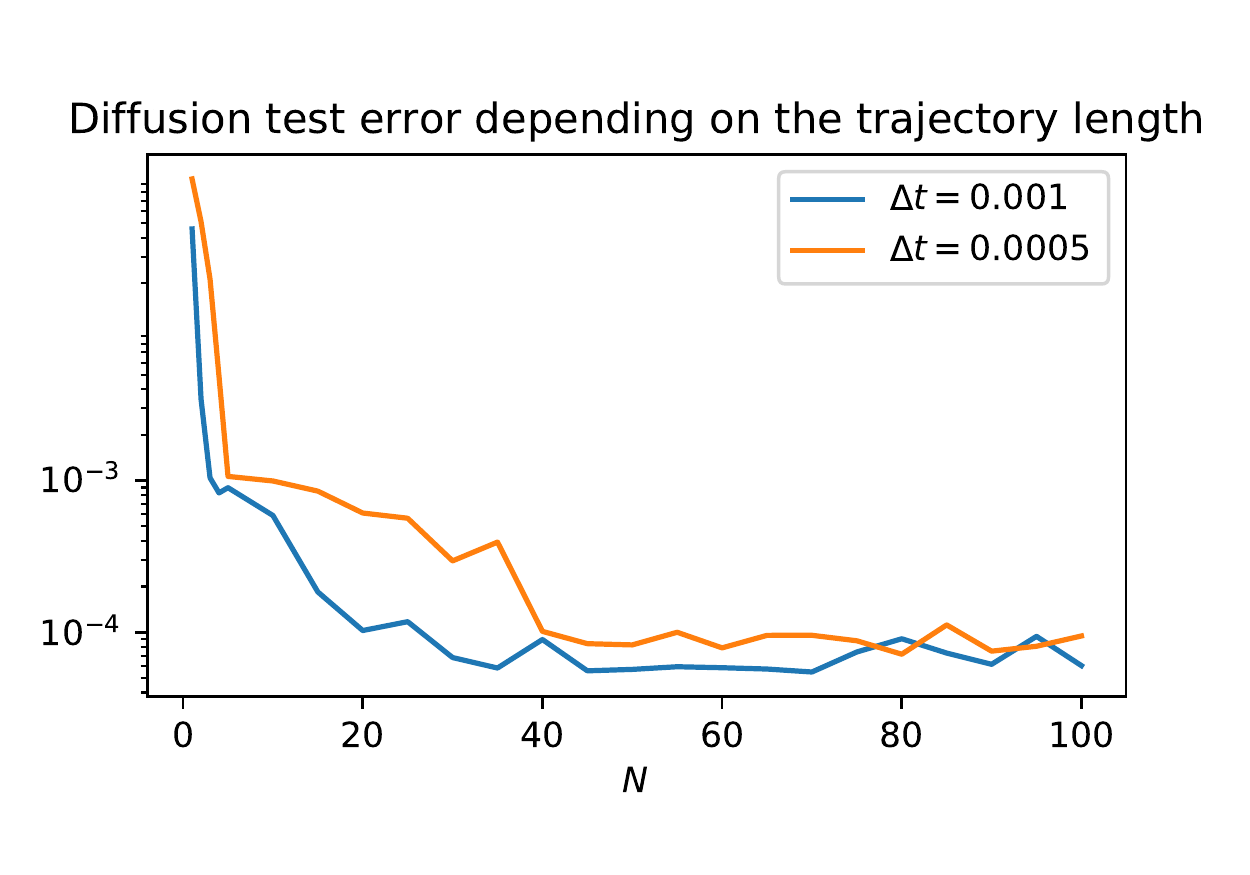}
  \caption{Weight $\alpha_\text{int} = 0.1$.}
  \label{fig: trajectory length diffusion loss weight 0.1}
\end{subfigure}%
\begin{subfigure}{.5\textwidth}
  \centering
  \includegraphics[width=\linewidth]{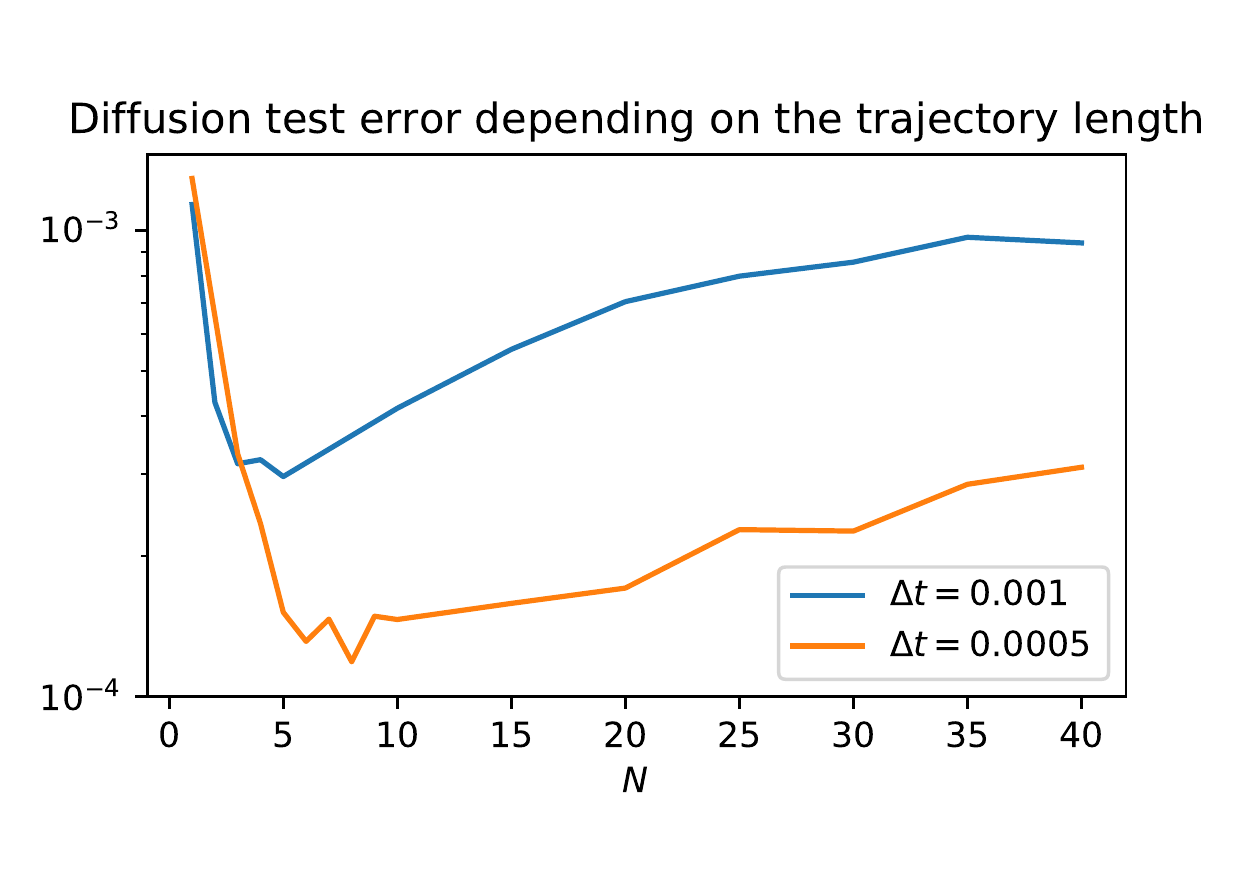}
  \caption{Weight $\alpha_\text{int} = 10$.}
  \label{fig: trajectory length diffusion loss weight 10}
\end{subfigure}
\caption{We display the $L^2$ error that one attains when using different choices of the maximal Euler steps $N$ in the diffusion loss for different discretization step-sizes $\Delta t$.}
\label{fig: trajectory length diffusion loss}
\end{figure}

In the diffusion loss as stated in \Cref{def: diffusion loss} we are free to choose the length $\mathfrak{t}$ of the forward trajectories, which affects the generated training data. Let us therefore investigate how different choices of $\mathfrak{t}$ influence the performance of \Cref{algorithm}. To this end, we consider again the elliptic problem from \Cref{sec: Elliptic problem with Dirichlet boundary data} and vary $\mathfrak{t}$. To be precise, let us fix different step-sizes $\Delta t$ and vary the Euler steps $N$ (recalling that $\mathfrak{t} = N \Delta t$), once choosing the weight $\alpha_{\text{int}} = 0.1$, as before, and once by considering $\alpha_{\text{int}} = 10$. For the former choice we can see in \Cref{fig: trajectory length diffusion loss weight 0.1} that larger trajectories tend to be better until a plateau is reached, whereas for the latter it turns out that there seems to be an optimal choice of the trajectory length, as displayed in \Cref{fig: trajectory length diffusion loss weight 10}. \par\bigskip

\subsubsection{Elliptic problem requiring the full Hessian matrix}
We consider the setting specified in \eqref{eq:toy}, replacing however the diffusion coefficient and the nonlinearity by
\begin{equation}
    \sigma(x, t) = \sqrt{\frac{2}{d}} \begin{pmatrix}
  1  & \cdots & 1 \\
  \vdots &  \ddots & \vdots \\
  1 & \cdots & 1
 \end{pmatrix}, \qquad h(x, y, z) = -2\gamma y \left(\gamma \sum_{i, j = 1}^d x_i x_j + d\right) + \sin\left(e^{2\gamma |x|^2} - y^2\right),
\end{equation}
respectively. Again, we can check that $V(x) = e^{\gamma |x|^2}$ is the unique solution to the corresponding boundary value problem. Since $\sigma$ is not diagonal anymore, the differential operator \eqref{eq: infinitesimal generator variational} contains the full Hessian matrix of second-order derivatives of the candidate solution $\varphi$. As discussed in \Cref{sec: from losses to algorithms} (see \Cref{tab: loss comparisons challenges}) this particularly impacts the runtime of the PINN method since all derivatives need to be computed explicitly. For the BSDE and diffusion losses, on the other hand, second-order derivatives are implicitly approximated using the underlying Brownian motion and we therefore do not expect significantly longer runtimes. \par\bigskip

\begin{figure}[H]
\centering
  \centering
  \includegraphics[width=1.0\linewidth]{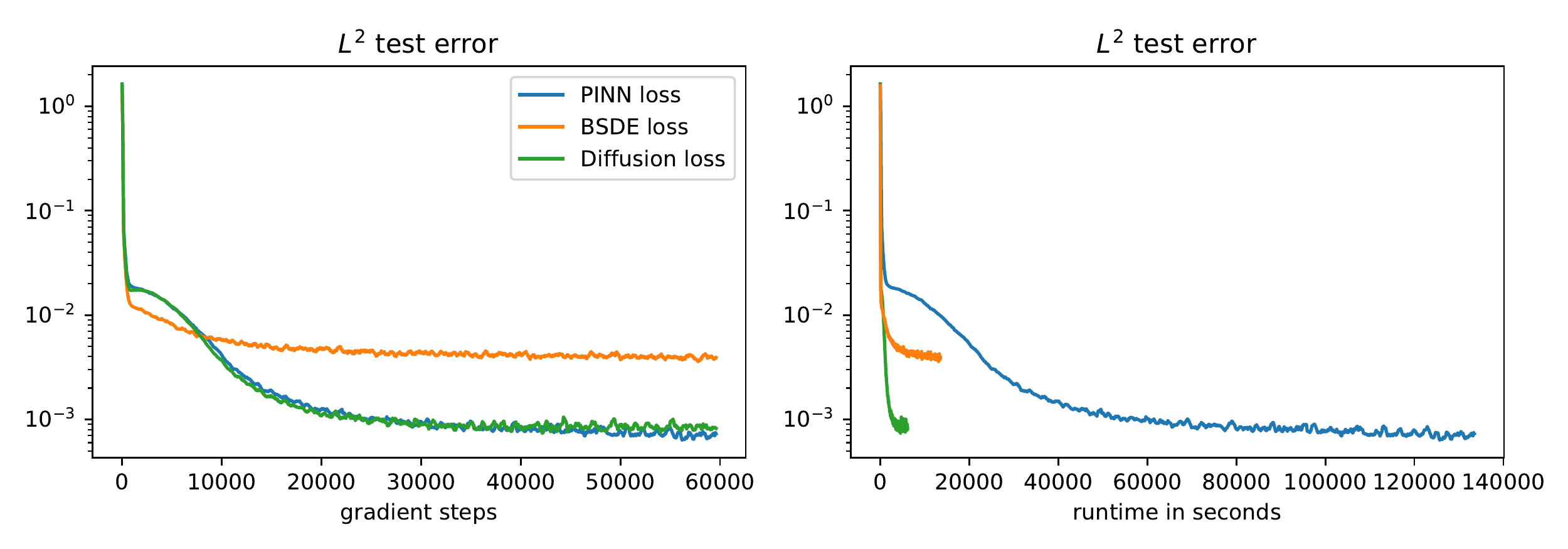}
  \caption{ $L^2$ error during the training process evaluated on test data for the three losses, plotted against the number of gradient steps (left panel) and  against the runtime (right panel). The diffusion loss relies on trajectory lengths of $\mathfrak{t} = 0.02$.}
  \label{fig: exponential on ball full Hessian}
\end{figure}

Let us consider $d = 20$ and $\gamma = 1$ and as before choose $N = 20$ (i.e. $\mathfrak{t} = 0.02$) for the diffusion loss. In \Cref{fig: exponential on ball full Hessian} we display the $L^2$ error during the training process, plotted against the number of gradient steps (left panel) and against the runtime (right panel). As expected, the PINN loss takes significantly longer. This effect should become even more pronounced with growing state space dimension $d$. 

\subsubsection{Parabolic problem with Neumann boundary data}

Let us now consider the parabolic problem \eqref{eq: general boundary value problem}, with the Dirichlet boundary condition \eqref{eq: space boundary condition} replaced by its Neumann counterpart
\begin{equation}
    \partial_{\vec{n}} V(x, t)  = g^N(x, t), \qquad   (x, t) \in \partial \Omega \times [0, T].
\end{equation} 
Here, $\partial_{\vec{n}} V:= \nabla V \cdot \vec{n}$ refers to the (outward facing) normal derivative at the boundary $\partial \Omega$.
We take
\begin{subequations}
  \begin{equation}
            b(x, t) = \mathbf{0}, \qquad \sigma(x, t) = \sqrt{2}\,  \mathrm{Id}_{d \times d}, \qquad f(x) = e^{\gamma |x|^2 + T}, \qquad g^N(x, t) = 2 \gamma e^{\gamma + t},\\
  \end{equation}
    \begin{equation}
      h(x, t, y, z) = -y (2\gamma  (2\gamma |x|^2 + d) + 1) + \sin\left(e^{2\gamma |x|^2 + 2t} - y^2\right).
  \end{equation}
  \end{subequations}
In this case,
\begin{equation}
    V(x, t) = e^{\gamma |x|^2 + t}
\end{equation}
provides the unique solution.\par\bigskip

\begin{figure}[H]
\centering
  \centering
  \includegraphics[width=1.0\linewidth]{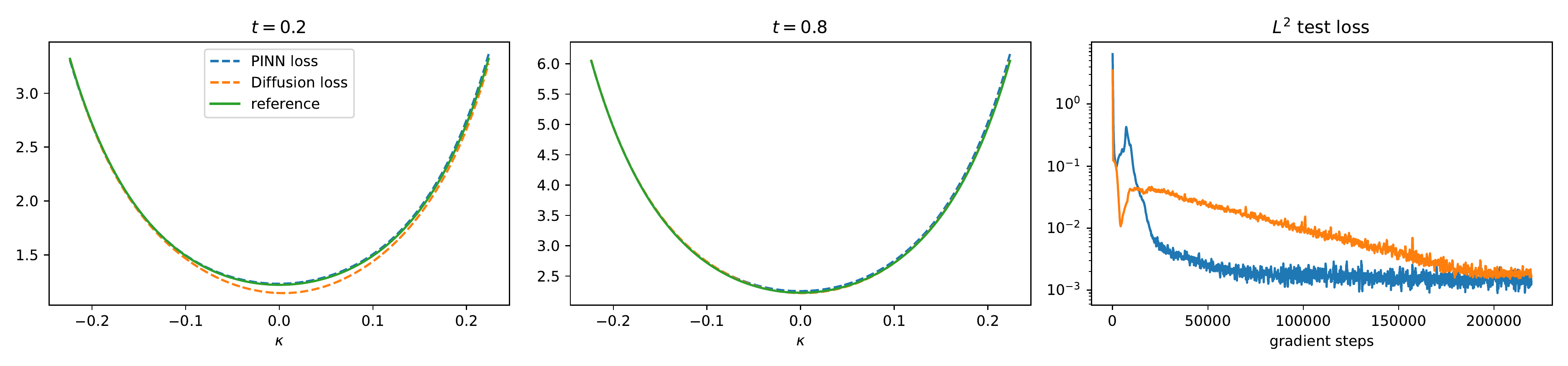}
  \caption{Left and central panel: Approximations along a curve for two different times using the diffusion and PINN losses. Right: $L^2$ test error along the training iterations.}
  \label{fig: exponential on ball parabolic Neumann}
\end{figure}

We choose $d = 20$ and $\gamma = 1$. For the diffusion loss we take $N = 25$ (i.e. $\mathfrak{t} = 0.025$) for the trajectory length. In the left and central panels of \Cref{fig: exponential on ball parabolic Neumann} we display the approximated solutions along the curve $\left\{(\kappa, \dots, \kappa)^\top : \kappa \in [-1 /\sqrt{d}, 1 /\sqrt{d} ] \right\}$ for two different times. We can see that both the diffusion and the PINN loss work well, with small advantages for the PINN loss. The right panel displays the $L^2$ test error over the iterations and confirms this observation.

\subsection{Committor functions}

Committor functions are important quantities in molecular dynamics as they specify likely transition pathways as well as  transition rates between (potentially metastable) regions or conformations of interest \cite{lu2015reactive,weinan2006towards}. Since for most practical applications those functions are high-dimensional and hard to compute, there have been recent attempts to approach this problem using neural networks \cite{khoo2019solving, li2019computing, rotskoff2020learning}. Based on the fact that committor functions fulfill elliptic boundary value problems, we can rely on the methods discussed in this paper.  \par\bigskip

For an $\mathbb{R}^d$-valued stochastic process $(X_t)_{t \ge 0}$ with continuous sample paths and two disjoint open sets $A,B \subset \mathbb{R}^d$, the committor function $V$ computes the probability of $X$ hitting $A$ before $B$ when starting in $x \in \mathbb{R}^d$, that is, 
\begin{equation}
    V(x) = \P(\tau_B < \tau_A | X_0 = x) = \E[\mathbbm{1}_B(X_\tau)|X_0 = x].
\end{equation}
Here, $\tau_A = \inf\{t > 0: X_t \in A \}$ and $\tau_B = \inf\{t > 0: X_t \in B \}$ refer to the hitting times corresponding to the sets $A$ and $B$. 
In the case when $(X_t)_{t \ge 0}$ is given as the unique strong solution to the SDE \eqref{eq: uncontrolled SDE}, it can be shown 
via the Kolmogorov backward PDE \cite[Section 2.5]{pavliotis2014stochastic} that $V$ fulfills the elliptic boundary value problem
\begin{align}
    LV = 0, \qquad& V\vert_{\partial A} = 0,\qquad V|_{\partial B} = 1,
\end{align}
where $L$ as in \eqref{eq: infinitesimal generator variational} refers to the associated infinitesimal generator,
see, for instance, \cite{weinan2006towards}. In the notation of \eqref{eqn: elliptic PDE} we have $\Omega = \R^d \setminus (A \cup B)$, $h = 0$ and $g(x) = \mathbbm{1}_B(x)$. \par\bigskip

\begin{figure}[H]
\centering
\includegraphics[width=0.8\linewidth]{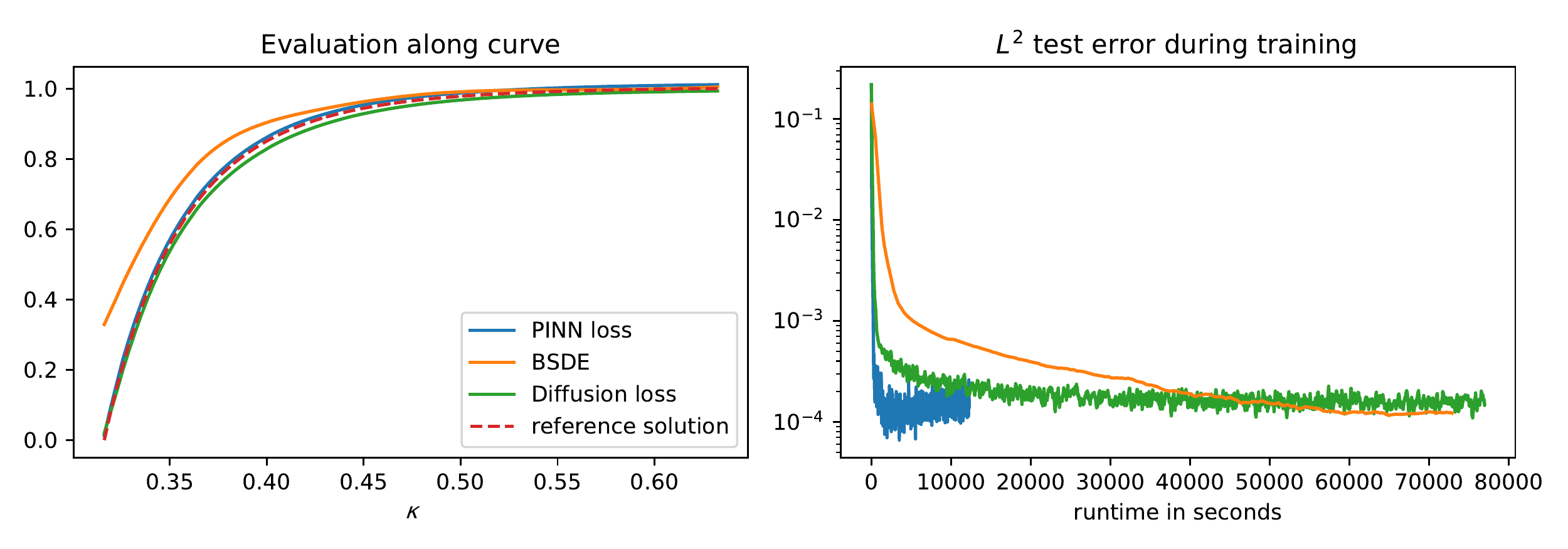}
\caption{Left: approximations of the $10$-dimensional committor function evaluated along a curve. Right: moving average of the test $L^2$ error along the training iterations.}
\label{fig: committor curve and loss}
\end{figure}

\begin{figure}[H]
\centering
\includegraphics[width=\linewidth]{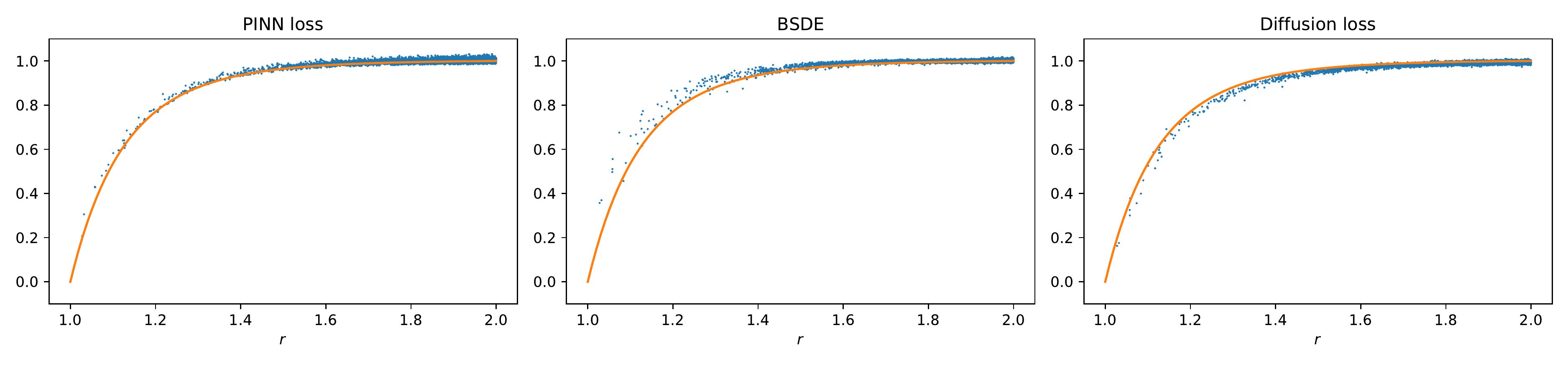}
\caption{We plot the approximated committor functions evaluated at $10000$ points uniformly sampled from the domain $\Omega$ (blue dots) and compare those to the reference solution (orange line) as a function of $r = |x|$.}
\label{fig: committor sampled evaluations}
\end{figure}

Following \cite[Section V.A]{hartmann2019variational}, we consider $(X_t)_{t \ge 0}$ to be a standard Brownian motion starting at $x \in \mathbb{R}^d$, that is,  $X_t = x + W_t$, corresponding to $b = \mathbf{0}$ and $\sigma = \mathrm{Id}_{d \times d}$ in \eqref{eq: uncontrolled SDE}. The sets $A$ and $B$ are defined as
\begin{align}
    A = \{x \in {\R}^d: |x| < a \},\qquad B = \{x \in {\R}^d:  |x| > b \},
\end{align}
with $b > a > 0$.
Hence, in this case the committor function describes the statistics of leaving a spherical shell through one of its boundaries. 
The solution takes the form
\begin{equation}
\label{eq: true committor solution}
    V(x) = \frac{a^2 - |x|^{2-d}a^2}{a^2-b^{2-d}a^2},
\end{equation}
for $d \ge 3$. Let us consider $d = 10$ as well as $a = 1$, $b = 2$. We take a DenseNet with the $\tanh$ activation function and compare the three losses against each other. For the diffusion loss we take $N = 50$ (i.e. $\mathfrak{t} = 0.05$) for the trajectory length. In \Cref{fig: committor curve and loss} we display the approximated solutions along a curve  $\left\{(\kappa, \dots, \kappa)^\top : \kappa \in [a /\sqrt{d}, b /\sqrt{d} ] \right\}$ in the left panel, realizing that in particular the PINN and diffusion losses lead to good approximations. This can also be observed in the right panel, where we plot a moving average of the $L^2$ error on test data based on a moving window of length $200$. The BSDE loss appears to be especially error-prone close to the left end-point of the curve displayed in \Cref{fig: committor curve and loss}, that is, close to the inner shell. Due to the volume distortion in high dimensions, few samples are drawn according to $\nu_\Omega$ for small values of $|x|$, see \Cref{fig: committor sampled evaluations}. We hence conclude that in this example, the BSDE loss suffers particularly from the relative sparsity of the samples, possibly in conjunction with numerical errors made while estimating the hitting times at the inner shell (see Section \ref{sec: SDE simulations}).

\subsection{Parabolic Allen-Cahn equation on an unbounded domain}

The Allen-Cahn equation in $d = 100$ has been suggested as a benchmark problem in \cite{weinan2017deep}. It is an example of a parabolic PDE posed on an unbounded domain,
\begin{subequations}
\begin{align}
    (\partial_t + L) V(x, t) + V(x, t) - V^3(x, t) &= 0,\qquad  & (x, t) \in \R^d \times [0, T], \\
    V(x, T) & = f(x), \qquad  & x \in \R^d,
\end{align}
\end{subequations}
with $f(x) = \left(2 + \frac{2}{5}|x|^2 \right)^{-1}$ and $T = \frac{3}{10}$.  We restrict attention to the ball $\{ x \in \R^d : |x| < r \}$ with radius $r = 7$, from which we sample the initial points of the trajectories for the BSDE and diffusion losses, or the training data for the PINN loss, respectively. Instead of using the uniform distribution on this set, we consider sampling uniformly on a box around the origin with side length $2$ and multiplying each data point $x$ by $\frac{r}{|x|}$. In contrast to uniform sampling, this approach generates more samples close to the origin, which we observe to slightly improve 
the accuracy of the obtained solutions. For the diffusion loss we choose $N = 25$ (i.e. $\mathfrak{t} = 0.025$) for the trajectory length.
We compare our approximations to a reference solution at $x_0 = (0, \dots, 0)^\top$ for different times $t \in [0, T]$ that is provided by a branching diffusion method specified in \cite{weinan2017deep}. In \Cref{fig: allen-cahn solutions} we see that all three attempts match this reference solution, with very minor advantages (for instance at the right end point) for the PINN and diffusion losses. In \Cref{tbl: computation time allen-cahn} we display the computation times until convergence, realizing that the BSDE loss needs significantly longer, which is in accordance with \Cref{rem:comparison BSDE and PINN} and might be explained by the fact that $\mathfrak{t} < T$. We note that the computation times are longer in comparison to e.g. \cite{weinan2017deep} since we aim for a solution on a given domain, whereas other attempts only strive for approximating the solution at a single point. \par\bigskip

\begin{minipage}{\textwidth}
  \begin{minipage}{0.49\textwidth}
    \centering
    \includegraphics[width=1.\linewidth]{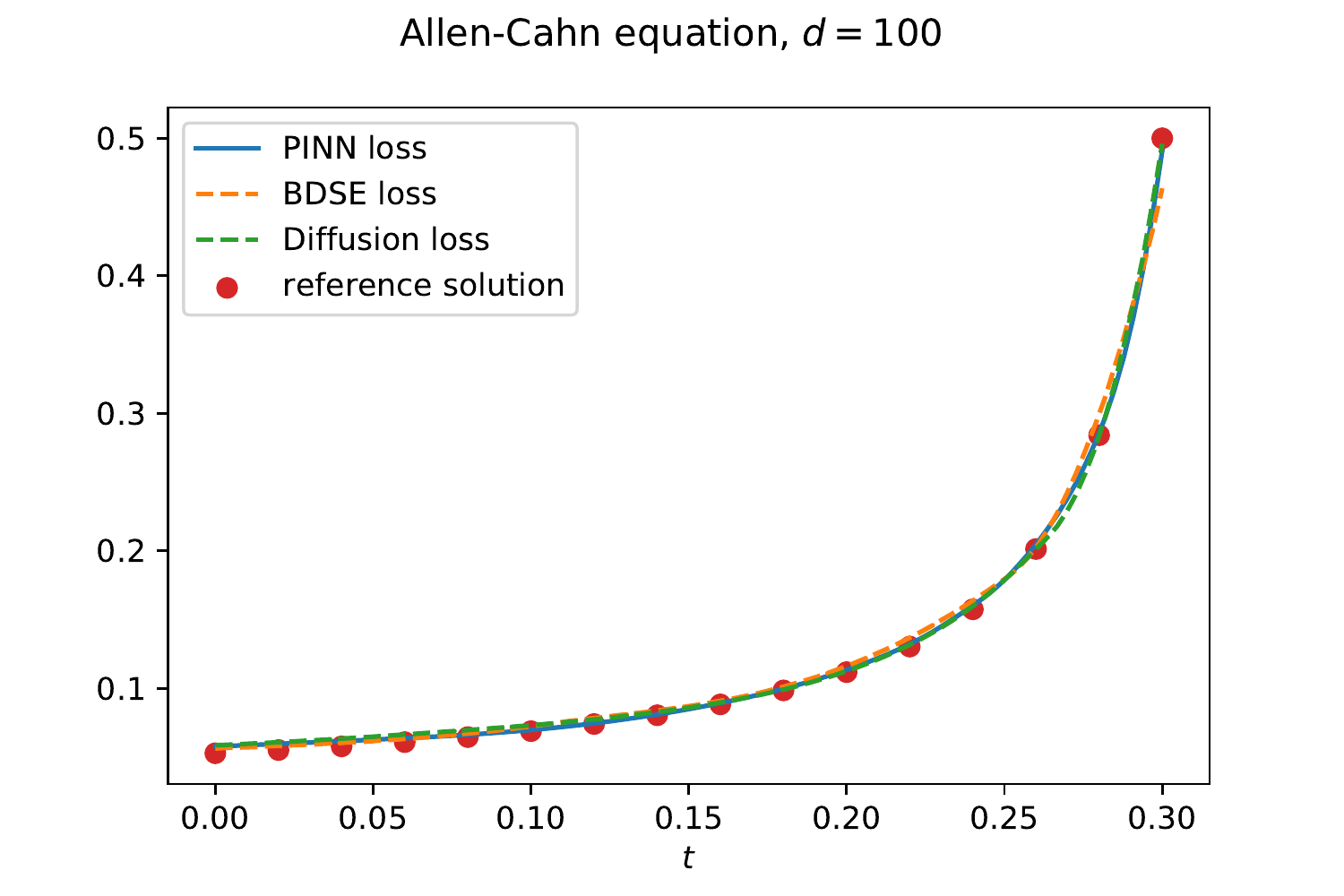}
    \captionof{figure}{Approximation of the solution to an Allen-Cahn equation in $d=100$ using different losses compared to a reference solution at $x_0 = (0, \dots, 0)^\top$ for different times $t \in [0, T]$.}
    \label{fig: allen-cahn solutions}
  \end{minipage}
  \hfill
  \begin{minipage}{0.49\textwidth}
    \centering
    \bgroup
    \def\arraystretch{1.3}
    \begin{tabular}{lr}
       \multicolumn{2}{c}{\textbf{Computation time}}  \\  \hline
        PINN loss & $325.46$ min\\  \hline
        BSDE loss & $4280.68$ min \\  \hline
        Diffusion loss & $194.38$ min
      \end{tabular}
      \egroup
      \captionof{table}{Computation times until convergence.}
      \label{tbl: computation time allen-cahn}
    \end{minipage}
  \end{minipage}

\subsection{Elliptic eigenvalue problems}

In this section we provide two examples for the approximation of principal eigenvalues and corresponding eigenfunctions. The first one is a linear problem and therefore \Cref{prop: eigenvalue loss} assures that the minimization of an appropriate loss as in \eqref{eq: eigenfunction loss} leads to the desired solution. The second example is a nonlinear eigenvalue problem, for which we can numerically show that our algorithm still provides the correct solution.

\subsubsection{Fokker-Planck equation}

As suggested in \cite{han2020solving}, we aim at computing the principal eigenpair associated to a Fokker-Planck operator,  defined by
\begin{equation}
\label{eq:FKP}
    LV = 
    -\Delta V - \nabla \cdot(V \nabla \Psi),
\end{equation}
for $V:\Omega \to \R$
on the domain $\Omega = [0, 2 \pi]^d$, and where $\Psi(x) = \sin\left(\sum_{i=1}^d c_i \cos(x_i) \right)$ is a potential with constants $c_i \in [0.1, 1]$, assuming periodic boundary conditions. This results in the eigenvalue problem
\begin{equation}
    \Delta V(x) + \nabla \Psi(x) \cdot \nabla V(x) + \Delta \Psi(x) V(x) = -\lambda V(x),
\end{equation}
 and 
 \begin{equation}
      V(x) = e^{-\Psi(x)}
 \end{equation}
is an eigenfunction associated to the principal eigenvalue $\lambda = 0$, see \cite[Section 4.7]{pavliotis2014stochastic}. We choose $c_i = 0.1$, $i =1,\ldots, d$, and approach this problem in dimension $d = 5$ following \Cref{sec: eigenvalue problems}, i.e. by minimizing the loss \eqref{eq: eigenfunction loss}, where for $\mathcal{L}$ we choose the diffusion loss with $N = 20$ (i.e. $\mathfrak{t} = 0.02$) and the periodic boundary condition is encoded via the term \eqref{eq: periodic boundary loss term}. Here and in the following eigenvalue problem the positivity of the approximating function is enforced by adding a ReLU function after the last layer of the DenseNet. \par\bigskip 

\begin{figure}[H]
\centering
  \centering
  \includegraphics[width=1.0\linewidth]{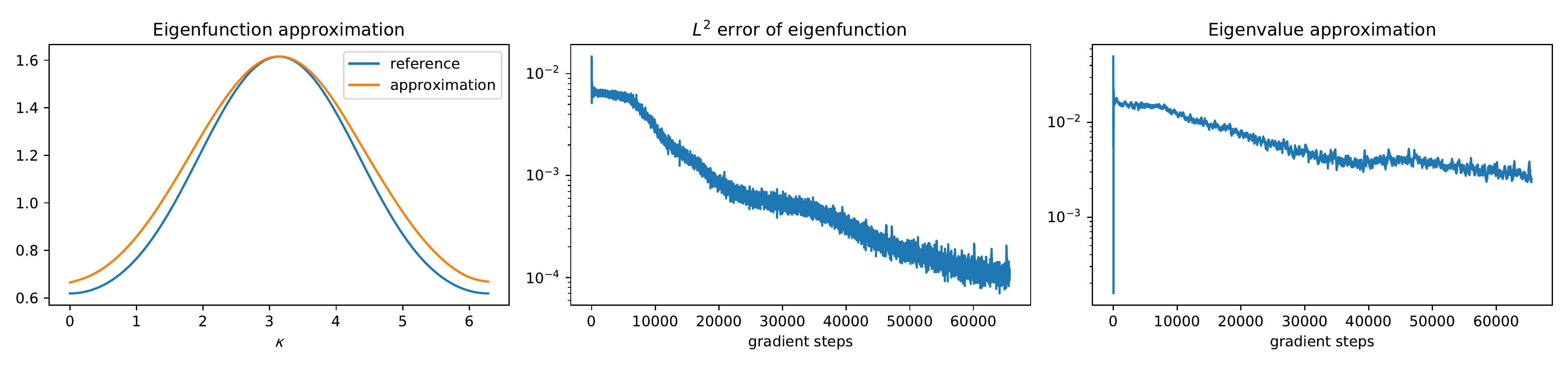}
  \caption{Left: Approximation and reference of the eigenfunction corresponding to the principal eigenvalue of the Fokker-Planck operator \eqref{eq:FKP}. Middle: $L^2$ error w.r.t. test data along the training iterations. Right: Moving average of the absolute value of the approximated eigenvalue along the gradient steps.}
  \label{fig: Fokker-Planck eigenpair d5}
\end{figure}

In the left panel of \Cref{fig: Fokker-Planck eigenpair d5} we display the approximated eigenfunction along the curve $\left\{(\kappa, \dots, \kappa)^\top : \kappa \in [0, 2\pi] \right\}$ and compare it to the reference solution. In the central panel we show the $L^2$ error w.r.t. the reference solution evaluated on uniformly sampled test data along the training iterations. The right panel displays the moving average of the absolute value of the eigenvalue taken over a moving window of $100$ gradient steps (since the true value is $\lambda = 0$ it is not possible to compute a relative error here). We see that both the eigenfunction and the eigenvalue are approximated sufficiently well.

\subsubsection{Nonlinear Schr\"odinger equation}
\label{sec:nonlinear sch}
\begin{figure}[H]
\centering
  \centering
  \includegraphics[width=1.0\linewidth]{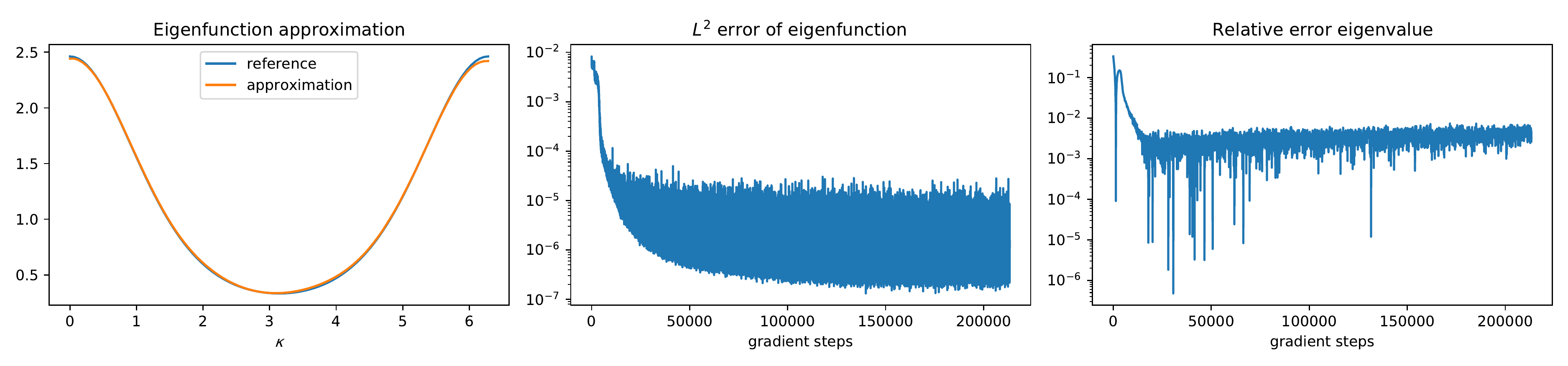}
  \caption{Left: Approximation and reference of the eigenfunction corresponding to the principal eigenvalue of the nonlinear Schr\"odinger operator in $d = 5$. Middle: $L^2$ error w.r.t. test data along the training iterations. Right: Relative error of the approximated eigenvalue along the gradient steps.}
  \label{fig: Schroedinger eigenpair d5}
\end{figure}

Let us now consider a nonlinear eigenvalue problem. Again following \cite{han2020solving}, we consider the nonlinear Schr\"odinger operator including a cubic term that arises from the Gross-Pitaevskii equation for the single-particle wave function in a Bose-Einstein condensate \cite{gross1961structure,pitaevskii1961vortex}. To be precise, we consider
\begin{equation}
    \Delta V(x) - V^3(x) - \Psi(x) V(x) = -\lambda V(x),
\end{equation}
where 
\begin{equation}
    \Psi(x) = -\frac{1}{c^2} \exp\left(\frac{2}{d} \sum_{i=1}^d \cos x_i \right) + \sum_{i=1}^d \left(\frac{\sin^2(x_i)}{d^2} - \frac{\cos x_i}{d} \right) - 3.
\end{equation}
One can show that 
\begin{equation}
    V(x) = \frac{1}{c} \exp \left(\frac{1}{d}\sum_{i=1}^d \cos x_i \right)
\end{equation}
is the eigenfunction corresponding to the principal eigenvalue $\lambda = -3$, where $c$ is determined from the normalization $\int_\Omega V^2(x) \,\mathrm dx = |\Omega|$. We add the latter constraint into the loss function by replacing the term $\mathcal{L}_c(\varphi)$ in \eqref{eq: eigenfunction loss} with $\mathcal{L}_n(\varphi) = \left(\E\left[\varphi(X)^2\right] - 1\right)^2$, where the distribution of $X$ is uniform on $\Omega = [0,2\pi]^d$. We again rely on the diffusion loss with a trajectory length of $N = 20$ (i.e. $\mathfrak{t} = 0.02$). \Cref{fig: Schroedinger eigenpair d5} shows the approximate solution of the eigenfunction in $d=5$ evaluated along the curve $\left\{(\kappa, \dots, \kappa)^\top : \kappa \in [0, 2\pi] \right\}$ as well as its $L^2$ error and the relative error of the approximated eigenvalue along the training iterations. We observe that both the eigenfunction and the eigenvalue are approximated quite well.

\begin{figure}[H]
\centering
  \centering
  \includegraphics[width=1.0\linewidth]{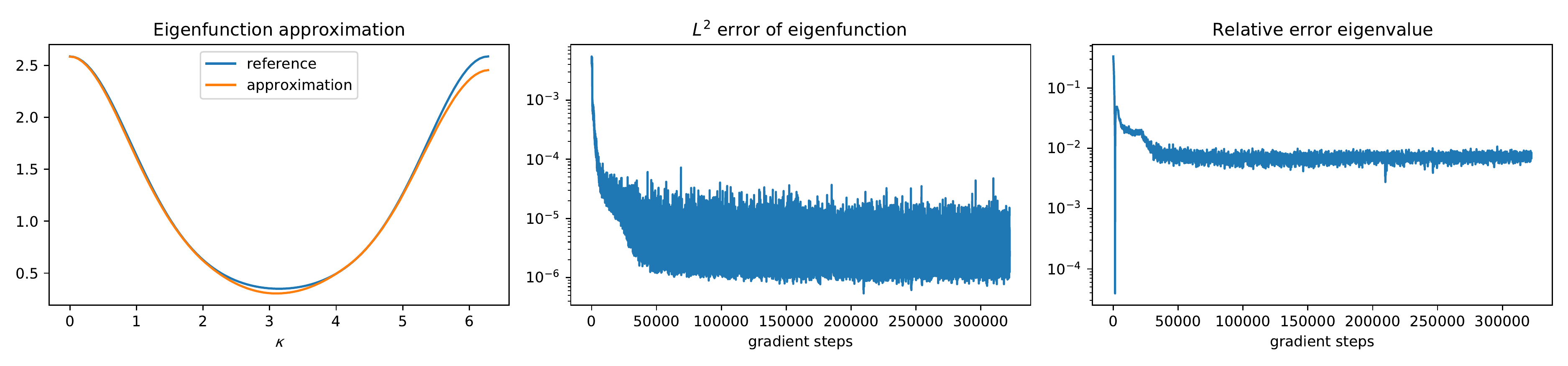}
  \caption{Same experiment as in \Cref{fig: Schroedinger eigenpair d5} in dimension $d=10$.}
  \label{fig: Schroedinger eigenpair d10}
\end{figure}

We repeat the experiment in dimension $d=10$ and display the results in \Cref{fig: Schroedinger eigenpair d10}.  The optimization task gets slightly more difficult, but the resulting eigenfunction and eigenvalue still fit the reference solution reasonably well. Note that in both experiments no explicit boundary conditions, but only the norm constraint and the periodicity are imposed.

\section{Conclusion and Outlook}
\label{sec: conclusion}

In this paper, we have investigated the relationship between BSDE and PINN based approximation schemes for high-dimensional PDEs through the lens of the novel diffusion loss. In particular, we have shown that the diffusion loss provides an interpolation between the aforementioned methods and demonstrated its promising numerical performance in a range of experiments, allowing us to trade-off strengths and weaknesses of BSDEs and PINNs. Although we believe that the diffusion loss may be a stepping stone towards a unified understanding of computational approaches for high-dimensional PDEs, many questions remain open: First, there is a need for principled insights into the mechanisms with which BSDE, PINN and diffusion losses can or cannot overcome the curse of dimensionality. Second, it is of great practical interest to optimize the algorithmic details, preferably according to well-understood theoretical foundations. In this regard, we mention
sophisticated (possibly adaptive) choices of the weights $\alpha_{\text{int}}$, $\alpha_{\text{b}}$, $\alpha_{\text{T}}$ and the measures $\nu_{\Omega \times [0,T]}$, $\nu_{\Omega}$, $\nu_{\partial \Omega \times [0,T]}$ as well as variance reduction techniques for estimator versions of the losses (see, for instance, \cite[Remark 4.7]{nusken2021solving}). Last but not least, we expect that the concepts explored in this paper may be fruitfully extended to the setting of parameter-dependent PDEs.
\\\\

\textbf{Acknowledgements.} This research has been funded by Deutsche Forschungsgemeinschaft (DFG) through the
grant CRC 1114 ‘Scaling Cascades in Complex Systems’ (projects A02 and A05, project number 235221301). We would like to thank Carsten Hartmann for many useful discussions.

\printbibliography

\end{document}